\documentclass[11pt]{article}
\usepackage{amsmath,amsfonts,amssymb,mathrsfs,cases,mathdots,colordvi,url,extarrows,graphicx}

\usepackage{relsize}
\usepackage{mathrsfs}
\usepackage{color}
\usepackage{indentfirst}
\usepackage{float}

\newtheorem{defn}{Definition}[section]
\newtheorem{assump}{Assumption}[section]
\newtheorem{lemma}{Lemma}[section]
\newtheorem{thm}{Theorem}[section]
\newtheorem{prop}{Proposition}[section]
\newtheorem{cor}{Corollary}[section]
\newtheorem{example}{Example}[section]
\newtheorem{remark}{Remark}[section]

\renewcommand{\Box}{\rule{2.2mm}{2.2mm}}

\def\beginproof{\par\noindent {\bf Proof.}\ \ }
\def\endproof{\hskip .5cm $\Box$ \vskip .5cm}

\def\beginproof{\par\noindent {\bf Proof.}\ \ }
\def\endproof{\hskip .5cm $\Box$ \vskip .5cm}

\begin{document}
\title{Optimality conditions for bilevel programs via Moreau envelope reformulation\thanks{The alphabetical order of the paper indicates the equal contribution to the paper.}}
\author{Kuang Bai\thanks{Department of Applied Mathematics, The Hong Kong Polytechnic University, Hong Kong, People's Republic of China. The research of this author was partially supported by the NSFC Grant 12201531 and the Hong Kong Research Grants Council PolyU153036/22p. Email: kuang.bai@polyu.edu.hk. } \and Jane J. Ye \thanks{Corresponding author. Department of Mathematics and Statistics, University of Victoria, Canada. The research of this author was partially
supported by NSERC. Email: janeye@uvic.ca.} \and Shangzhi Zeng\thanks{Department of Mathematics and Statistics, University of Victoria, Canada. Email: zengshangzhi@uvic.ca.}}
  \date{}
\maketitle 
\begin{abstract} 
	For bilevel programs with a convex lower-level program, the classical approach  replaces the lower-level program with its Karush-Kuhn-Tucker  condition and solve the resulting mathematical program with complementarity constraint (MPCC). It is known that when the set of lower-level  multipliers is not unique,  MPCC may not be equivalent to the original bilevel problem, and many MPCC-tailored constraint qualifications do not hold. In this paper, we study bilevel programs where the lower level is generalized convex. Applying the equivalent reformulation via Moreau envelope, we derive new directional optimality  conditions.  Even in the nondirectional case, the new optimality condition is stronger than the strong stationarity  for the corresponding MPCC.

	\vskip 10 true pt
	

	\vskip 10 true pt
	
	
\end{abstract}

\section{Introduction}

In this paper, we consider the standard bilevel program:
\begin{align*}
\mbox{(BP)}\quad \min_{x,y} \quad F(x,y)\
s.t.\ \ y\in S(x),\ G(x,y)\leq0,
\end{align*}
where for any given $x$, $S(x)$ denotes the solution set of the lower-level program
\[
(\mathrm{P_x})\quad \min_{y} f(x,y) \ \ \mbox{s.t. }  g(x,y)\leq 0.
\] 
In this paper unless otherwise specified we assume that $F, f:\mathbb R^n\times\mathbb R^m\rightarrow\mathbb R,\ G:\mathbb R^n\times\mathbb R^m\rightarrow\mathbb R^q$, $ g:\mathbb R^n\times\mathbb R^m\rightarrow\mathbb R^p$ are continuously differentiable. Throughout this paper, by saying that   the lower-level program is generalized convex we mean that the following assumption holds.
\begin{assump}[Generalized convexity of the lower-level program] For each $x$, $f$ is convex in $y$ and $g$ is quasi-convex in variable $y$, i.e., the level set
$\left \{y| g(x,y)\leq r \right \}$ is convex for all real number $r$.
\end{assump}
 Under the lower-level generalized convexity if a certain constraint qualification for (P$_x$) holds at $y$, then  $y\in S(x)$ if any only if  there {\it exists} a Lagrange multiplier $\lambda$ satisfying the Karush-Kuhn-Tucker (KKT) condition, i.e.,
\[
 \nabla_y f(x,y) +\nabla_y g(x,y)^T \lambda=0,\ 0\leq\lambda\perp g(x,y)\leq 0.
\] Denote the set of all Lagrange multipliers of (P$_x$) at $y$ by $\Lambda(x,y)$. If the lower-level program is generalized convex and a certain constraint qualification holds, then (BP) is equivalent to the following single-level problem
\begin{eqnarray*}
\quad	\begin{array}{ll}
		\min_{x,y}~ &F(x,y)\\
		s.t.~ &G(x,y)\leq0,
		\mbox{ there exists } \lambda \in \Lambda(x,y).
	\end{array}	
\end{eqnarray*}
 However the above problem is not track-able.  Hence instead of considering the above equivalent problem, it became popular to study the following related problem
\begin{eqnarray*}
	{\rm (MPCC)}\quad
	\begin{array}{ll}
		\min_{x,y,\lambda}~ &F(x,y)\\
		s.t.~ &G(x,y)\leq0,\
		 \nabla_y f(x,y) +\nabla g(x,y)^T \lambda=0,\
		0\leq\lambda\perp g(x,y)\leq 0,
	\end{array}	
\end{eqnarray*} which is a mathematical program with complementarity constraint (MPCC).
This approach is the so-called  first-order approach or the KKT approach, which was popularly studied over the last three decades; see e.g. \cite{Dempe,Luo,Outrata} for the general theory  and  \cite{DZ2012,GY19,Ye00,Ye05,YY} for the optimality conditions derived by using this approach. 
Now in addition, we assume the functions $f,g$ are twice continuously differentiable and recall the S-stationarity condition for (MPCC):
	\begin{defn}\label{SMPCC}\cite[Definition 2.7]{Ye05}
	A vector $(\bar x,\bar y,\bar\lambda)$  is said to be strong (S-)stationary of (MPCC) if there exists $(\mu_G,\mu_g,\mu_{e},\mu_\lambda)\in\mathbb R^{q+p+m+q}$ such that
	\begin{equation}\label{S-cond}
	\begin{aligned}
		&\nabla F(\bar x,\bar y)+\nabla G(\bar x,\bar y)^T\mu_G+\nabla g(\bar x,\bar y)^T\mu_g+\nabla_{x,y}  l(\bar x,\bar y,\bar \lambda)^T\mu_{e}=0,\\
		&\nabla_\lambda l(\bar x,\bar y,\bar \lambda)^T\mu_{e}-\mu_\lambda=0,\\
        &0\leq\mu_G\perp G(\bar x,\bar y)\leq 0,\ \mu_\lambda \perp\bar\lambda,\
        \mu_g\perp g(\bar x,\bar y)\leq 0,\\
       &  \nabla_y f(\bar x,\bar y) +\nabla g(\bar x,\bar y)^T \bar \lambda=0,\
		0\leq \bar \lambda\perp g(\bar x,\bar y)\leq 0,\\
        &(\mu_g)_j\geq0\ \mbox{and}\ (\mu_\lambda)_j\geq0,\ \forall j \in  I_{0}(\bar x,\bar y):=\{j|\bar\lambda_j=0, g_j(\bar x,\bar y)=0\},
	\end{aligned}
	\end{equation}
	where 
	$l(x,y,\lambda):=\nabla_y f(x,y) +\nabla_y g(x,y)^T \lambda.$
\end{defn}
It is well-known that  if $(\bar x,\bar y,\bar\lambda)$ is  a local optimal solution of (MPCC) satisfying  the MPCC linear independent constraint qualification (MPCC-LICQ) \cite[Definition 2.8]{Ye05}, then it must be S-stationary. 

Although (MPCC) is more track-able than the original bilevel program, there are two major drawbacks of using the KKT approach. First, when the set of multipliers $\Lambda(x,y)$ is not a singleton for certain feasible solution $(x,y)$,  (MPCC) is usually not equivalent to (BP) if local optimality is considered.  Specifically, if  $(\bar x,\bar y,\bar \lambda)$ is a local minimizer of (MPCC),   $(\bar x,\bar y)$ may not be a local minimizer of (BP); see the details in Dempe and Dutta \cite{DD}. Secondly, MPCC variant of Mangasarian-Fromovitz constraint qualification (MPCC-MFCQ) may not hold at feasible points of (MPCC).
In fact for the case where $g(x,y)$ is independent of $x$, Gfrerer and Ye \cite{GY17} showed that when $\Lambda(x,y)$ is not a singleton for all feasible $(x,y)$, {MPCC-MFCQ} is never satisfied by any feasible solution of (MPCC). Furthermore, in this case even the  MPCC Guignard constraint qualification (MPCC-GCQ), one of the weakest MPCC-tailored constraint qualifications may not be valid.

Another way to reformulate (BP) as a single-level problem is the value function approach, which was first proposed in \cite{VP} for numerical purpose.
Using the optimal value function  $v(x):=\inf_{y}\{f(x,y)| g(x,y)\leq0\}$ of the lower-level program, the constraint $y\in S(x)$ can always be equivalently replaced by inequalities $f(x,y)-v(x)\leq 0, g(x,y)\leq 0$. Ye and Zhu \cite{YZ95} studied the optimality condition of the resulting single-level program
\begin{align*}
	\mbox{(VP)}~~~~\quad &\min_{x,y} \quad F(x,y)\\
	&s.t.\quad f(x,y)-v(x)\leq 0,\ g(x,y)\leq 0,\ G(x,y)\leq0.
\end{align*}
Under the Lipschitz continuity of the value function and a certain constraint qualification, at a feasible solution $(\bar x,\bar y)$, the nonsmooth {stationarity} condition in terms of Clarke subdifferential can be written as follows: there exist $(\alpha, \lambda_g, \lambda_G)\in \mathbb R_+\times \mathbb R^{p+q}$ such that
\begin{equation}
	\begin{aligned}\label{KKTVP}
	&0\in  \nabla F(\bar x,\bar y)+\alpha( \nabla f(\bar x,\bar y)-\partial^c v(\bar x)\times\{0\})+\nabla g(\bar x,\bar y)^T\lambda_g+\nabla G(\bar x,\bar y)^T\lambda_G,\\
	  &0\leq \lambda_g\perp g(\bar x,\bar y)\leq 0,\ 0\leq \lambda_G\perp G(\bar x,\bar y)\leq 0,
\end{aligned} 
\end{equation} where  $\partial^c v$ denotes the Clarke subdifferential  of the value function $v$. However, the value function is implicitly defined, and to develop a verifiable optimality condition, there are two major challenges. First, it is known that KKT conditions may not hold at a local minimizer unless a certain constraint qualification holds. Since the strict inequality $f(x,y)-v(x)<0$ is nowhere satisfied in the feasible region, the no nonzero abnormal multiplier constraint qualification (NNAMCQ) fails at any feasible point of (VP). This means  constraint qualifications weaker than NNAMCQ need to be employed. Secondly, $\partial^c v$ is implicitly defined.
Under appropriate conditions, say the restricted inf-compactness condition (see e.g. \cite[Definition 3.8]{GLYZ}) plus the MFCQ at each optimal point, one has the Lipschitz continuity of the value function and an upper estimate can be stated as follows:
$$\partial^c v(\bar x) \subseteq \operatorname{co} \bigcup_{y\in S(\bar x)} \bigcup_{\lambda\in \Lambda  (\bar x, y)} \{\nabla_x f(\bar x, y)+\nabla_x g(\bar x, y)^T\lambda \},$$  
where {$\operatorname{co} C$} denotes the convex hull of set $C$.
In the literature (see e.g., \cite{YZ95}), to get an explicit optimality condition, one would  replace the set $\partial^c v(\bar x)$ by its upper estimates. 
Hence, to verify the optimality condition, one has to check each point in $S(\bar x)$. Moreover, the resulting {stationarity} system is not computable by any algorithms if $S(\bar x)$ is  not a singleton. A natural question is: in application, can we find all elements in $S(\bar x)$? The answer is yes when bilevel polynomial programs are considered, see e.g., \cite{NieYe17, NieYe21}. But for general (BP), it is hard to calculate the set $S(\bar x)$. A natural idea is to replace $S(\bar x)$ by a smaller set. If  either the functions $f,g$ are convex in both variables $x$ and $y$ or the inner semicontinuity (see e.g., \cite{Aub2}) of $S(x)$ at $(\bar x,\bar y)$ holds, then  $S(\bar x)$ in the union is replaced by $\{\bar y\}$, i.e.,
\begin{align*}\partial^c v(\bar x) \subseteq  \bigcup_{\lambda\in \Lambda  (\bar x, \bar y)} \{\nabla_x f(\bar x, \bar y)+\nabla_x g(\bar x, \bar y)^T\lambda \}.\end{align*}
 By replacing $\partial^c v(\bar x)$ in (\ref{KKTVP}) with an element in the above upper estimate,  
 the following KKT {stationarity} condition can be obtained: there exist $(\alpha, \lambda_g, \lambda_G)\in \mathbb R_+\times \mathbb{R}^{p+q}$ and $\bar \lambda\in \Lambda(\bar x,\bar y)$ such that
\begin{equation}\label{eqn2}
	\begin{aligned}
	&0=\nabla  F(\bar x,\bar y)+\nabla g(\bar x,\bar y)^T (\lambda_g -\alpha  \bar \lambda)+\nabla G(\bar x,\bar y)^T\lambda_G,\\
&0\leq \lambda_g\perp g(\bar x,\bar y)\leq 0,\ 0\leq \lambda_G\perp G(\bar x,\bar y)\leq 0.
	\end{aligned}
	\end{equation}
	Let $ \mu_g:=\lambda_g -\alpha \bar \lambda.$ Then 
	$(\mu_g)_j\geq 0$ if $\bar\lambda_j=0$. 
Hence $(\mu_G,\mu_g,\mu_e,\mu_\lambda):=(\lambda_G,\lambda_g -\alpha \bar \lambda,0,0)$ satisfies the {stationarity} condition (\ref{S-cond}). Therefore the {stationarity} condition (\ref{eqn2}) is stronger than the S-{stationarity} condition (\ref{S-cond}) for problem (MPCC). Moreover unlike in the MPCC reformulation, the multiplier $\bar  \lambda$ is not arbitrarily given, it is    selected from the set of multipliers $\Lambda(\bar x,\bar y)$. {The optimality condition in the form (\ref{eqn2}) was first derived in \cite[Theorem 4.2]{Ye04} (see also \cite[Theorem 4.1 and Corollary 4.1]{Ye06} and \cite[Theorem 8.8.3]{Yebook}, \cite[Lemma 5.1(c)]{LM2023b}) under the joint convexity assumption on the lower-level program (i.e., $f, g$ convex in both variables $x,y$)},  or under the semicontinuity of the solution set (see e.g. \cite[Theorem 6.21]{Boris}) and  is  used in some recent work on semismooth Newton method for solving bilevel programs \cite{FisherZemkohoZhou,LM2023b,LM2023a}. 

However, the joint convexity  of the lower-level program or the inner semicontinuity property of the solution map is quite strictive. So the problem is left to be solved: are there other conditions that also allow us to derive the KKT condition as the simple system (\ref{eqn2}).

Recently, Gao et al. \cite{DwC2023} proposed a modified value function reformulation for (BP).
For any given $x,y, \gamma>0$, we define the {\it Moreau envelope} function
of the lower-level program as
	\begin{equation*}
		\begin{aligned}
			v_{\gamma}(x,y) & :
			= \inf_{w } \left\{\left. f(x,w) + \frac{1}{2\gamma}
			\|w - y\|^2 \right| g(x,w) \le 0  \right\}.
		\end{aligned}
	\end{equation*}
	When the lower-level program $(\mathrm{P_x})$ is generalized convex, the constraint $y\in S(x)$ can be equivalently replaced by $f(x,y)-v_\gamma(x,y)\leq0, g(x,y)\leq0$. Hence, (BP) can be equivalently reformulated as:
	\begin{equation*}
		{\rm (VP)}_\gamma~~~~~~~~~~\begin{aligned}
			\min_{x,y}  ~~& F(x,y) \\
			s.t. ~~~~& f(x,y)-v_\gamma (x,y)\leq 0,\ g(x,y)\leq0, \
G(x,y)\leq0.
		\end{aligned}
	\end{equation*}
Provided that $v_\gamma$ is Lipschitz continuous, one can immediately write out the KKT {stationarity} condition at a feasible solution $(\bar x,\bar y)$ for (VP)$_\gamma$:  there exist $(\alpha, \lambda_g, \lambda_G)\in \mathbb R_+\times \mathbb{R}^{p+q}$ such that
\begin{equation}\label{KKTMR}
\begin{aligned}
	&0\in \nabla F(\bar x,\bar y) +\alpha (\nabla f(\bar x,\bar y)-\partial^c v_\gamma(\bar x,\bar y))+\nabla g(\bar x,\bar y)^T\lambda_g+\nabla G(\bar x,\bar y)^T\lambda_G,\\
	&0\leq \lambda_g\perp g(\bar x,\bar y)\leq 0 ,\ 0\leq \lambda_G\perp G(\bar x,\bar y)\leq 0.
\end{aligned} 
\end{equation}
Generally, it is easier to study KKT condition (\ref{KKTMR})  than  (\ref{KKTVP}). Indeed, when $\gamma$ is sufficiently small, provided that the function $f(x,w)$ is weakly convex for each fixed $x$, the objective function $f(x,w) + \frac{1}{2\gamma}
\|w - y\|^2$ is strongly convex in $w$. It follows that for each $x,y$, the objective function for the minimization problem in the definition of $v_\gamma(x,y)$ is coercive and the solution map is always a singleton provided that $g$ is quasi-convex in $w$. As a result, many regularity conditions for sensitivity analysis of $v_\gamma$ hold automatically, and under relatively weak assumptions, say the MFCQ of the system $g(\bar x,y)\leq 0$ at $\bar y\in S(\bar x)$, we have the sharper upper estimate:
$$ \partial^c v_\gamma(\bar x,\bar y) \subseteq \bigcup_{\lambda\in \Lambda  (\bar x,\bar y)} \left \{(\nabla_x f(\bar x,\bar y)+\nabla_x g(\bar x,\bar y)^T\lambda,0) \right \}.$$
It is interesting that replacing $\partial^c v_\gamma(\bar x,\bar y)$ by the above upper estimate still results in  the  KKT {stationarity} system (\ref{eqn2}).  Besides, the above upper estimate for $v_\gamma$ is sharper than that for the value function $v(x)$ in the absence of full convexity; see e.g., \cite{BY2023,Gauvin}. Hence, the resulting optimality condition is stronger than classical optimality conditions for (VP). 

Another direction to improve the value function approach is to consider the directional constraint qualifications and KKT conditions. Recently, Bai and Ye \cite{BY} derived directional KKT condition for (VP) under the directional calmness condition, where the value function $v(x)$ is assumed to be directionally Lipschitz continuous and directionally differentiable along a certain direction, and the solution set $S(x)$ is replaced by the solution along certain critical direction which is in general is a smaller subset of $S(x)$. Since these directional conditions are usually weaker than their non-directional counterparts, the resulting directional KKT condition is in general sharper and the directional constraint qualification is weaker than the nondirectional counterpart. 

In this paper, by using the Moreau envelope reformulation, we derive a directional version of the KKT system (\ref{eqn2})  under the directional calmness condition. We also derive some verifiable directional quasi-normality condition as sufficient conditions for the directional calmness condition.  Note that the directional version of the KKT system is in general stronger than the nondirectional KKT system while the directional calmness condition is in general weaker than the nondirectional calmness condition. A direct consequence of our result gives some new conditions under which the KKT system (\ref{eqn2}) holds without assuming either full convexity of the lower-level problem or the inner semicontinuity of the solution map $S(x)$.

We organize the paper as follows. In the next section, we provide notations, preliminaries and preliminary results. In Section 3 we derive the directional KKT condition under the directional calmness condition for a nonsmooth nonlinear optimization problem. In Section 4, we study directional sensitivity analysis of the Moreau envelope function.  In Section 5, we apply the results of previous sections to (VP)$_\gamma$ and derive a verifiable necessary optimality condition. We conclude the paper in Section 6.


\section{Preliminaries and preliminary results} We first give notations that will be used in the paper. We denote by $\overline{\mathbb R}:=\mathbb R\cup\{\pm\infty\}$, while $\mathbb R:=(-\infty,+\infty)$. $\|\cdot\|$ denotes the Euclidean norm. $\langle a,b\rangle$ denotes the inner product of vectors $a,b$, and $a\perp b$ means the inner product $\langle a,b\rangle=0$. Let $\Omega$ be a set. By $x^k\xrightarrow{\Omega}\bar{x}$ we mean $x^k\rightarrow\bar{x}$ and for each $k$, $x^k\in \Omega$.  By $x^k\xrightarrow{u}\bar x$ where $u$ is a vector, we mean that the sequence $\{x^k\}$ approaches $\bar x$ in direction $u$, i.e., there exist $t_k\downarrow 0, u^k\rightarrow u$ such that $x^k=\bar x+ t_k u^k$. By {$\operatorname{o}(t)$}, we mean $\lim_{t\rightarrow0}\frac{o(t)}{t}=0$.
We denote by $\mathbb B$, $\bar{\mathbb B}$, $\mathbb S$  the open unit ball, the closed unit ball and 
  the unit sphere, respectively. $\mathbb B_\delta(\bar z)$ denotes the open ball centered at $\bar z$ with radius $\delta$.
We denote by {$\operatorname{ri}\Omega$, $\operatorname{int}\Omega$, $\operatorname{co}\Omega$ and $\operatorname{cl}\Omega$} the relative interior, the interior,  the convex hull and  the closure of a set $\Omega$, respectively. The distance from a point $x$ to a set $\Omega$ is denoted by  ${\operatorname{dist}(x,\Omega):=\inf\{\|x-y\||y\in\Omega\}}$, the indicator function of the set $\Omega$ is denoted by $\delta_\Omega(x)$, and the support function of set $\Omega$  is denoted by  $\sigma_\Omega(x):=\sup\{\langle x,y \rangle|y\in\Omega\}$. For a single-valued map $\phi:\mathbb R^n\rightarrow\mathbb R^m$, we denote by $\nabla \phi(x)\in \mathbb{R}^{m\times n}$   the Jacobian matrix of $\phi$  at $x$ and for a function $\varphi:\mathbb R^n\rightarrow\mathbb R$, we denote by $\nabla \varphi(x)$ the gradient of $\varphi$  at $x$. Denote the pre-image of set $\Omega\subseteq\mathbb R^m$ under map $\phi$ by $\phi^{-1}(\Omega):=\{x\in\mathbb R^n|\phi(x)\in\Omega\}$. For a set-valued map $\Phi:\mathbb R^n\rightrightarrows\mathbb R^m$ the graph of $\Phi$ is defined by  ${\operatorname{gph}\Phi:=\{(x,y)| y\in \Phi(x)\}}$. For an extended-valued function $\varphi:\mathbb R^n\rightarrow \overline{\mathbb R}$, we define its domain by ${\operatorname{dom}\varphi:=\{x\in\mathbb R^n| \varphi(x) <\infty\}}$, and its epigraph by ${\operatorname{epi}\varphi:=\{(x,\alpha)\in\mathbb R^{n+1}|\alpha\geq\varphi(x)\}}$. For a function $g:\mathbb R^n\rightarrow\mathbb R$, we denote $g_+(x):=\max \{0,g(x)\}$ and if it is vector-valued then the maximum is taken componentwise.


In the following, we gather some basic concepts and prove some results in variational analysis, which will be used later on. For more details see e.g. \cite{BHA,BS,Clarke,DonRock,Long,Aub2,RW}.
\begin{defn}[Tangent Cone and Normal Cone] (see, e.g., \cite[Definitions 6.1 and 6.3]{RW})
	Given a set $\Omega\subseteq\mathbb{R}^n$ and a point $\bar{x}\in \Omega$, the tangent cone to $\Omega$ at $\bar{x}$ is defined as
	$$T_\Omega(\bar{x}):=\left \{d\in\mathbb{R}^n|\exists t_k\downarrow0, d_k\rightarrow d\ \mbox{ s.t. } \bar{x}+t_kd_k\in \Omega\ \forall k\right \}.$$
	The regular normal cone and the limiting normal cone  to $\Omega$ at $\bar{x}$ are defined as
\begin{eqnarray*}\widehat{N}_\Omega(\bar{x})&:=&\left \{ \zeta\in \mathbb{R}^n\bigg| \langle \zeta ,x-\bar{x}\rangle \leq o(\|x-\bar x\|) \quad \forall x\in \Omega \right \},\\
	 N_\Omega(\bar{x})&:=&\left \{\zeta\in \mathbb{R}^n\bigg| \exists \ x_k\xrightarrow{\Omega}\bar{x},\ \zeta_k{\rightarrow}\zeta\ \text{such that}\ \zeta_k\in\widehat{N}_\Omega(x_k)\ \forall k\right \},\end{eqnarray*}
	respectively.
	\end{defn}

	\begin{defn}[Directional Normal Cone] (\cite[Definition 2.3]{GM}  or \cite[Definition 2]{Gfr13}). 
    Given a set $\Omega \subseteq \mathbb{R}^n$, a point $\bar x \in \Omega$ and a direction $ u \in\mathbb{R}^n$, the limiting normal cone to $\Omega$ at $\bar{x}$ in direction $u$ is defined by
    $$N_\Omega(\bar{x};u):=\left \{\zeta \in \mathbb{R}^n\bigg| \exists \ t_k\downarrow0, u_k\rightarrow u, \zeta_k\rightarrow\zeta  \mbox{ s.t. } \zeta_k\in \widehat{N}_\Omega(\bar{x}+t_k u_k)\ \forall k \right \}.$$
\end{defn}
It is obvious that $N_{\Omega}(\bar x; 0)=N_{\Omega}(\bar x)$, $N_{\Omega}(\bar x; u)=\emptyset$ if $u \not \in T_\Omega(\bar x)$ and $N_{\Omega}(\bar x;u)\subseteq N_\Omega(\bar x)$.
Moreover when $\Omega$ is convex, by \cite[Lemma 2.1]{Gfr14} the directional and the classical normal cone have the following relationship
\begin{equation}
N_\Omega(\bar x;u)=N_\Omega(\bar x)\cap \{u\}^\perp  \qquad \forall u\in T_\Omega(\bar x).\label{convNormal}
\end{equation}

When $u=0$ the following  definition coincides with the Painlev\'e-Kuratowski inner/lower and outer/upper  limit of $\Phi$ as $x\rightarrow\bar x$ respectively; see e.g., \cite{Aub2}.
\begin{defn}
	Given a set-valued map $\Phi:\mathbb R^n\rightrightarrows\mathbb R^m$ and a direction $
u \in\mathbb R^n$, the inner/lower and outer/upper limit of $\Phi$ as $x\xrightarrow{u}\bar x$ respectively is defined by
	\begin{align*}
	\liminf_{x\xrightarrow{u}\bar x} \Phi(x) &:= \{y\in\mathbb R^m| \forall \ \mbox{sequence}\ t_k\downarrow0, u^k\rightarrow u, \exists y^k\rightarrow y \mbox{ s.t. } y^k\in \Phi(\bar x+t_k u^k)\}, \\
	\limsup_{x\xrightarrow{u}\bar x} \Phi(x) &:= \{y\in\mathbb R^m| \exists   \ \mbox{sequence}\ t_k\downarrow0, u^k\rightarrow u,  y^k\rightarrow y \mbox{ s.t. } y^k\in \Phi(\bar x+t_ku^k)\},
	\end{align*}
	respectively.
\end{defn}

\begin{defn}[Graphical and Directional Derivatives]\cite[Page 199]{RW}
Let $\phi:\mathbb R^n\rightarrow\mathbb{R}^m $ and $x, u\in \mathbb R^n$. The graphical derivative of $\phi$ at $x$ is 
\begin{equation*}
D \phi(x)(u):=\left \{v\left|\exists t_k\downarrow0, u_k\rightarrow u \ s.t.~ v=\lim_{k\rightarrow+\infty}\frac{\phi(x+t_ku_k)-\phi(x)}{t_k} \right.\right \}.
\end{equation*}
When $m=1$, the upper-/lower-Dini derivative of $\phi$ at $x$ in direction $u$ is
\begin{align*}
&\phi'_+(x;u):=\limsup_{t\downarrow0}\frac{\phi(x+tu)-\phi(x)}{t},\\
&\phi'_-(x;u):=\liminf_{t\downarrow0}\frac{\phi(x+tu)-\phi(x)}{t}.
\end{align*}
The usual directional derivative of $\phi$ at $x$ in the direction $u$ is 
$$\phi'(x;u):=\lim_{t\downarrow0}\frac{\phi(x+tu)-\phi(x)}{t}$$
when this limit exists. 
\end{defn}

Below we review the concept of directional neighborhood introduced by Gfrerer in \cite{Gfr13}.
Given a direction $u\in \mathbb{R}^n$, and positive numbers $\epsilon,\delta>0$, the directional neighborhood of direction $u$ is a set defined by 
\begin{equation*}
	{\cal V}_{\epsilon,\delta}(u):=\{x\in\epsilon\mathbb{B}|\big\| \| u\| x-\| x\| u\big\|\leq\delta \|x\| \|u\|\}.
\end{equation*}
Clearly, the directional neighborhood of direction $u=0$ is just the open ball $\epsilon\mathbb{B}$ and the directional neighborhood of a nonzero direction $u\not =0$ is a smaller subset of $\epsilon\mathbb{B}$.
Hence many regularity conditions can be extended to a directional version which is weaker than the original nondirectional one. 

We say that $\phi:\mathbb R^n\rightarrow\mathbb{R}^m $ is directionally Lipschitz continuous  at $\bar x$ in direction $u$ if there are positive numbers $L, \epsilon, \delta$ such that
$$\|\phi(x)-\phi(x') \|\leq L \|x-x'\|  \quad \forall x, x' \in \bar x+{\cal V}_{\epsilon,\delta}(u). $$
It is easy to see that if $\phi:\mathbb R^n\rightarrow\mathbb{R}^m $ is directionally Lipschitz continuous at $\bar x$ in direction $u$, then $D\phi(\bar x)(u)\neq\emptyset$ and both $\phi'_+(\bar x;u)$ and $\phi'_-(\bar x;u)$ are finite and belong to $D\phi(\bar x)(u)$.
Futhermore, if $\phi$ is directionally differentiable at $\bar x$ in direction $u$
then for any sequence $\{u^k\}$ which converges to $u$, we have
\begin{equation*}
\phi'(\bar x;u)=\lim_{k\rightarrow\infty}\frac{\phi(\bar x+t_ku^k)-\phi(\bar x)}{t_k}.
\end{equation*}

We now recall the definition of some subdifferentials below.
\begin{defn}[Subdifferentials](\cite[Definition 8.3]{RW})
Let  $\varphi:\mathbb{R}^n \rightarrow \overline{\mathbb{R}}$ and  $\varphi (\bar x)$ be finite. 
The Fr\'{e}chet (regular) subdifferential of $\varphi$ at $\bar{x}$ is the set
\begin{eqnarray*}
\widehat\partial \varphi(\bar{x}):=\left\{\xi\in\mathbb{R}^n|\varphi(x)\geq \varphi(\bar x)+\langle\xi,x-\bar x\rangle+{o(\|x-\bar x\|)}\right\},
\end{eqnarray*}
the limiting (Mordukhovich or basic) subdifferential of $\varphi$ at $\bar{x}$ is the set
\begin{align*}
\partial \varphi(\bar{x}):=
\{\xi\in\mathbb{R}^n|\exists x^k\rightarrow  \bar x,  \ \xi^k\rightarrow\xi \ \mbox{ s.t. }  \varphi(x^k)\rightarrow \varphi(\bar x), \xi^k\in \widehat\partial \varphi(x^k)\}.
\end{align*}
\end{defn}
Recently a directional version of the subdifferential was introduced \cite{Gfr13,GM} and further studies and related results can be found  in \cite{Long,BHA,BY2023}.
\begin{defn}[Directional Subdifferentials]\label{ads}
	  Let $\varphi:\mathbb{R}^n\rightarrow \overline{\mathbb{R}}$ and   $\varphi (\bar x)$ be finite. 
The limiting subdifferential of $\varphi$ at $\bar x$ in direction 
$u\in \mathbb{R}^n$ is defined as 
	\begin{eqnarray*}
	\partial \varphi(\bar x;u):=\{\xi\in\mathbb{R}^n|\exists t_k\downarrow0, u^k\rightarrow u, \xi^k\rightarrow\xi \mbox{ s.t. } \varphi(\bar x+t_ku^k)\rightarrow \varphi(\bar x), \   \xi^k\in \widehat{\partial} \varphi(\bar x+t_ku^k)\}.
	\end{eqnarray*} 
\end{defn} It is easy to see that if $u\notin T_{\operatorname{dom}\varphi}(\bar x)$, then $ \partial \varphi(\bar x;u)=\emptyset$ and $\partial \varphi(\bar x;0)=\partial \varphi(\bar x)$.   When $u\not =0$, $\partial \varphi(\bar x;u)\subseteq \partial \varphi(\bar x)$ and when $\varphi$ is continuously differentiable at $\bar x$,  $\partial \varphi(\bar x;u)=\{ \nabla \varphi(\bar x)\}$ for any direction $u$.

{ By \cite[Corollary 4.1]{BHA} directional subdifferentials can be equivalently defined by the directional normal cone to its epigraph. 
\begin{prop}\label{dsnorm}Let $\varphi:\mathbb{R}^n\rightarrow \overline{\mathbb{R}}$ and   $\varphi (\bar x)$ be finite.  Suppose $\varphi$ is Lipschitz continuous and directionally differentiable at $\bar x$ in direction $u$. Then
\begin{align*}
	\partial \varphi(\bar x;u)&=\{\xi\in\mathbb{R}^n|(\xi,-1)\in N_{\operatorname{epi} \varphi}(\bar x,\varphi(\bar x);(u,\varphi'(\bar x;u)))\}. \nonumber\\
	\end{align*}
\end{prop}
}
\begin{defn}[Directional Clarke Subdifferential]\cite[Definition 2.7]{BY}\label{Clarke}
Let $\varphi:\mathbb R^n\rightarrow\mathbb R$ be directionally Lipschitz continuous at  $\bar x$ in direction $u\in\mathbb R^n$. The directional Clarke subdifferential of $\varphi$ at $\bar x$ in direction $u$ is defined as
\begin{equation*}
\partial^c \varphi(\bar x;u):=\operatorname{co}(\partial \varphi(\bar x;u)).
\end{equation*}	
\end{defn}
 It is clear that  the directional Clarke subdifferential in direction $u=0$ coincides with the Clarke   subdifferential $\partial^c \varphi(\bar x)$. By \cite[Proposition 2.1]{BY}, one has 
 $\partial^c(-\varphi)(\bar x;u)=-\partial^c \varphi(\bar x;u).$

A function $\varphi(x):\mathbb{R}^n\rightarrow \overline{\mathbb{R}}$ is said to be weakly convex if there exists a modulus $\rho\geq 0$ such that the function $\varphi(x)+\frac{\rho}{2} \|x\|^2$ is convex. The class of weakly convex functions is large and it includes the class of functions with a bounded Lipschitz gradient, i.e., $\operatorname{dom}\varphi$ is open, $\varphi (x)$ is  differentiable at all $x \in \operatorname{ dom} \varphi$ and there exists $L>0$ such that $\|\nabla \varphi(x)-\nabla \varphi(y)\|\leq L\|x-y\|$ for $\forall x,y \in \operatorname{ dom} \varphi$. It turns out that like convex functions, weakly convex functions are directionally differentiable and Lipschitz continuous in the interior of its domain. Moreover the directional limiting subdifferential and the directional Clarke subdifferential coincide as shown below.
{
\begin{prop} \label{Prop2.2}  Let $\varphi:\mathbb{R}^n\rightarrow \overline{\mathbb{R}}$ be weakly convex and $\bar x\in \operatorname{ri}(\operatorname{dom}\varphi)$. Then 
one has
\begin{eqnarray}
&& \varphi'(\bar x;u)=\sup \left \{ \langle \xi, u\rangle |\xi \in \partial \varphi(\bar x)\right \}.\label{Cdirectionald}
\end{eqnarray}  Moreover if $\bar x\in \operatorname{int}(\operatorname{dom}\varphi)$, then $\varphi$ is directionally differentiable at $\bar x$ and Lipschitz continuous around $\bar x$, (\ref{Cdirectionald}) holds with ``sup'' replaced with ``max'' and 
\begin{eqnarray}
&& \partial^c \varphi(\bar x;u)= 
\partial \varphi(\bar x;u)=\partial\varphi(\bar x)\cap\{\xi|\xi^Tu=\varphi'(\bar x;u)\}. \label{C-L}\end{eqnarray}
\end{prop}
\beginproof By the weak convexity, $\tilde \varphi(x):=\varphi(x)+\rho\|x\|^2$ is convex { for some $\rho>0$}. Hence the limiting and regular subdifferential coincide with the one in the sense of convex analysis and we have
$$\partial \varphi(\bar x)=\partial \tilde \varphi(\bar x)-2\rho \bar x.$$ By \cite[Theorem 23.4]{rockafellar},  we have
$$  \varphi'(\bar x;u)=\tilde \varphi'(\bar x;u) -2\rho\bar x^Tu =\sup \left \{ \langle \xi, u\rangle |\xi \in \partial \tilde \varphi(\bar x)\right \}-2\rho\bar x^Tu$$ 
from which (\ref{Cdirectionald}) follows.
Now suppose $\bar x\in \operatorname{int}(\operatorname{dom}\varphi)$. Then by \cite[Theorem 23.4]{rockafellar}, $\tilde  \varphi( x) $ is directionally differentiable and hence $\varphi$ is directionally differentiable. And the convex function $ \tilde\varphi(x)$ is Lipschitz continuous around $\bar x$ and so is $\varphi(x)$. 
}

We now prove (\ref{C-L}). By sum rule of directional subdifferential \cite[Theorem 5.6]{Long},
	$$ \partial \varphi(\bar x;u)= \partial \tilde\varphi(\bar x;u) -2\rho \bar x.$$
By Proposition \ref{dsnorm}, $\partial  \tilde \varphi(\bar x;u)=\{\xi|(\xi,-1)\in N_{\operatorname{epi} \tilde \varphi}((\bar x,\tilde \varphi(\bar x));(u,\tilde \varphi'(\bar x;u))\}$. 
By the convexity of $\tilde \varphi$, $\operatorname{epi} \tilde \varphi$ is convex, so by (\ref{convNormal}) $$N_{\operatorname{epi}\tilde \varphi}((\bar x,\tilde\varphi(\bar x));(u,\tilde \varphi'(\bar x;u))=N_{\operatorname{epi}\tilde \varphi}((\bar x,\tilde \varphi(\bar x))\cap \{(u,\tilde \varphi'(\bar x;u))\}^\perp.$$
Then $\xi\in  \partial\tilde\varphi(\bar x;u)$ if and only if
	\[
	(\xi,-1)\in N_{\operatorname{epi}\tilde \varphi}((\bar x,\tilde \varphi(\bar x))\cap \{(u,\tilde \varphi'(\bar x;u))\}^\perp
	\] which means $\xi\in  \partial \tilde\varphi(\bar x;u)$ if and only if $\xi\in  \partial \tilde\varphi(\bar x)$ and $\xi^Tu=\tilde \varphi'(\bar x;u)$. Note that
	$$\tilde \varphi'(\bar x;u)=\varphi'(\bar x;u)+2\rho {\bar x}^T u.$$ Consequently, $\xi\in  \partial \tilde\varphi(\bar x;u)$ if and only if $\xi-2\rho \bar x \in  \partial \varphi(\bar x) $ and $(\xi-2\rho {\bar x})^Tu=\varphi'(\bar x;u) $. { This proves the second equality of (\ref{C-L}).}
	$\partial \varphi(\bar x)$ is a closed convex set, and so is $\partial \varphi(\bar x;u)$. By Definition \ref{Clarke}, $\partial^c \varphi(\bar x;u)=\partial \varphi(\bar x;u)$.
\endproof

We now give the definition of directional metric subregularity constraint qualification.
\begin{defn}[Directional MSCQ] \cite[Definition 2.1]{Gfr13}
Let $\bar z$ be a solution to the system $\phi(z)\leq0$, where $\phi:\mathbb R^a\rightarrow\mathbb R^b$.
 Given a direction $u\in\mathbb{R}^a$, the system $\phi(z)\leq0$ is said to satisfy the  metric subregularity constraint qualification (MSCQ) at $\bar{z}$ in direction $u$, or the set-valued map $\Phi(z):=\phi(z)-\mathbb R^b_-$ is metrically subregular at ($\bar z,0$)  in direction $u$, if  there are positive reals $\epsilon>0,\delta>0$ and $\kappa>0$ such that
		\begin{equation*}
		{\operatorname{dist}}(z,\phi^{-1}(\mathbb R^b_-))\leq\kappa \|\phi_+(z)\| \qquad \quad \forall  z\in\bar{z}+{\cal V}_{\epsilon,\delta}(u).
		\end{equation*}
\end{defn} 
If $u=0$ in the above definition, then we say that the system $\phi(z)\leq0$ satisfies MSCQ at $\bar z$. 

{It is known that MSCQ holds when MFCQ holds or $\phi(z)$ is affine.} In the following, we list some sufficient conditions for the directional MSCQ. The reader can refer to \cite{BY,BYZ,BHA,BCH} for more details.
\begin{defn}\label{qp} Let $\bar z$ be a solution to the system $\phi(z)\leq0$, where $\phi:\mathbb R^a\rightarrow\mathbb R^b$. Given a direction $u\in\mathbb{R}^a$.
	\begin{itemize}
		\item  Suppose that  $\phi$ is Lipschitz at $\bar z$. We say that the no-nonzero abnormal multiplier constraint qualification (NNAMCQ) holds at $\bar z$ if 
		\begin{equation*}
			0\in \partial\langle\zeta,\phi\rangle(\bar z) \mbox{ and } 0\leq\zeta\perp \phi(\bar z)  \Longrightarrow \zeta=0.
		\end{equation*}
		\item  Suppose that   $\phi$ is directionally Lipschitz at $\bar z$ in direction $u$. We say that the first order sufficient condition for metric subregularity (FOSCMS) holds  at $\bar z$ in direction $\xi$ where $\xi \in D\phi(\bar z)(u)$ if 
		there exists no $\zeta\neq 0$ satisfying 
	$$	0\in \partial\langle\zeta,\phi\rangle(\bar z;u) ,  0\leq\zeta\perp \phi(\bar z), \zeta\perp \xi.$$
		\item Suppose that  $\phi$ is directionally Lipschitz at $\bar z$ in direction $u$. We say that the directional quasi-normality  holds at $\bar{z}$ in direction $u$
		if 
		there exists no $\zeta\neq 0$ satisfying 
	$$	0\in \partial\langle\zeta,\phi\rangle(\bar z;u) ,  0\leq\zeta\perp \phi(\bar z), \zeta\perp \xi  \mbox{ for some } \xi\in D\phi(\bar z)(u) $$ and there exist sequences $t_k\downarrow0,\ u^k\rightarrow u$ satisfying
		\begin{equation*}
			\phi_j(\bar{z}+t_ku_k)>0,\ \mbox{ if }  j \in \{j=1,\ldots,b|\phi_j(\bar z)=0 \mbox{ and } \zeta_j>0\}.
		\end{equation*}
	\end{itemize}
\end{defn} 

It is known that in the smooth setting, NNAMCQ is equivalent to the MFCQ. { FOSCMS holding in a nonzero direction is in general weaker than MFCQ. }

{Throughout this paper, for any parametric inequality system $ \phi(x,y)\leq 0$ where $\phi(x,y): \mathbb{R}^n\times \mathbb{R}^m \rightarrow \mathbb{R}^s$ and $(\bar x,\bar y)\in \{(x,y) \mid \phi(x,y)\leq 0\}$,  we denote the active set of the inequality at $(\bar x,\bar y)$ as
	$$I_\phi(\bar x,\bar y):=\{i=1,\dots, s \mid \phi_i(\bar x,\bar y)=0\}.$$ }

The following regularity condition will be useful in this paper.
\begin{defn}[RCR Regularity](\cite[Definition 1]{MS}). 
	We say that the feasible  map ${\cal F}(x):=\{y|g(x,y)\leq0\}$  is relaxed constant rank (RCR) regular at $(\bar x,\bar y ) \in \operatorname{gph }{\cal F}$
	if there exists $\delta>0$ such that for any index subset $K\subseteq I_g(\bar x,\bar y)
	$, the family of gradient vectors $\nabla_y g_j(x,y), j\in K$, has the same rank at all points $(x,y) \in \mathbb{B}_\delta(\bar x,\bar y)$.
\end{defn}

Recall the directional Robinson Stability \cite[Definition 4.7]{BY}).
\begin{defn}[Directional Robinson Stability]\label{dRS}
	We say that the feasible map ${\cal F}(x)$ satisfies Robinson stability ${\rm (RS)}$ property at $(\bar x,\bar y)\in \operatorname{gph }{\cal F}$ in direction $u\in\mathbb R^n$ if there exist positive scalars $\kappa, \epsilon, \delta$ such that 
	\begin{equation*}
		\operatorname{dist}(y, {\cal F}(x))\leq \kappa \|g_+(x,y)\| \quad \forall x\in \bar x+{\cal V}_{\epsilon,\delta}(u), y\in  \mathbb{B}_\epsilon(\bar y).
	\end{equation*}
\end{defn}
If RS holds at $(\bar x,\bar y)$ in direction $u=0$, we say that  RS holds at $(\bar x,\bar y)$ (\cite[Definition 1.1]{HM}). 
RS means that the MSCQ holds at $\bar y$ uniformly in a directional neighborhood of $\bar x$ in direction $u$.

 Define the image directional derivative of $g$ with respect to $x$ at $(\bar x,\bar y)$ in direction $u$ as the closed cone $${\rm Im}D_xg(\bar x,\bar y;u):=\left\{\alpha v\left|\alpha\geq0, \exists t_k\downarrow0, u^k\rightarrow u\ s.t.\ v=\lim_{k\rightarrow\infty}\frac{g(\bar x+t_ku^k,\bar y)-g(\bar x,\bar y)}{\|g(\bar x+t_ku^k,\bar y)-g(\bar x,\bar y)\|}\right.\right\}.$$
We always have $$\{\alpha\nabla_xg(\bar x,\bar y)u|\alpha\geq0\}\subseteq {\rm Im}D_xg(\bar x,\bar y;u)$$ and the equality holds if $\nabla_xg(\bar x,\bar y)u\neq0$.

The following proposition lists some sufficient conditions for directional RS. The reader is referred to \cite{BY2023-1} for more sufficient conditions for directional RS. {And for any $\tilde u\in \mathbb{R}^n$, we need the following cone
 	\begin{equation} \label{Lconen} \mathbb L(\bar x,\bar y; \tilde  u):=\left \{v|\nabla g_i (\bar x,\bar y)( \tilde u,v)\leq0, \ i\in I_g(\bar x,\bar y) \right \}.\end{equation}}

\begin{prop}\label{Prop2.4new}\cite[Proposition 4.1]{BY}, \cite[Proposition 3.2, Corollaries 3.1 and 3.2]{BY2023-1}. 
	If the system $g(x,y)\leq0$ satisfies one of the following conditions at $(\bar x,\bar y)$, then RS holds at $(\bar x,\bar y)$ in direction  $u$.

	\begin{itemize}
		\item[1.] $g(x,y)=h(x)+By+c$, where $h:\mathbb R^n\rightarrow\mathbb R^p$ and $B$ is an {$p\times m$ matrix,}
			and the feasible region $\mathcal F(x)$ is nonempty for all $x\in\mathcal V_{\epsilon,\delta}(u)$ for some $\epsilon,\delta>0$.
			\item[2.] For any direction $d\in {\rm Im} D_xg(\bar x,\bar y; u)$, there exists $v\in \mathbb{R}^m$
satisfying 
\begin{equation}\label{eqnIm}
d_i+\nabla_y g_i(\bar x,\bar y)v\leq0,\ i\in I_g(\bar x,\bar y).\end{equation}
 FOSCMS for the system $g(\bar x, y)\leq0$ at $y=\bar y$ holds in direction $d+\nabla_y g(\bar x,\bar y)v$, i.e.,
\begin{eqnarray*}
	\left\{
	\begin{array}{ll}
		\nabla_yg(\bar x,\bar y)^T\lambda=0,\\
		0\leq\lambda\perp g(\bar x,\bar y),\ \lambda\perp (d+\nabla_y g(\bar x,\bar y)v )
	\end{array}
	\right.\implies \lambda=0,
\end{eqnarray*}
{for each $(d,v)\not =0$ satisfying (\ref{eqnIm}) and $d\in {\rm Im} D_xg(\bar x,\bar y; u)$.}
\item[3.] $\nabla_x g(\bar x,\bar y)u \not =0$. For any $\alpha\geq 0$, $ \mathbb L(\bar x,\bar y;\alpha  u)$ defined in (\ref{Lconen}) is nonempty and FOSCMS for the system $g(\bar x, y)\leq0$ at $y=\bar y$ holds { in all directions in $\mathbb L(\bar x,\bar y; \alpha u)$} for $\alpha \geq 0$.
	\item[4.]  MFCQ for the system $g(\bar x,y)\leq0$ holds at $y=\bar y$. 
	\end{itemize}
\end{prop}
\section{Directional KKT conditions}\label{Sec3}
In this section we derive directional KKT condition for  the optimization problem
\begin{eqnarray*}
{\rm (P)}~~~ \min_{z} && \varphi(z)\quad \mbox{ s.t. } \phi(z)\leq 0,
 \end{eqnarray*}
where $\varphi:\mathbb{R}^a \rightarrow \mathbb{R} $ and $\phi:\mathbb{R}^a \rightarrow \mathbb{R}^b $. In \cite[Section 3]{BY}, Bai and Ye obtained the directional optimality condition for (P) under the assumption that all constraint functions  $\phi_j$ are directionally Lipschitz continuous and directionally differentiable. In this section we extend the result to allow  at most one of the active constraint functions, say $\phi_1$ to be not directionally differentiable.

Let $\bar z$ be a feasible solution to problem (P). { Assume $1\in I_\phi(\bar z)$} throughout this section. Now suppose that  $\varphi$ is continuously differentiable at $\bar z$ and $\phi_j$ is directionally Lipschitz continuous for all $j$, and directionally differentiable for all $j\neq1$ at $\bar z$. Then  {treating the inequality constraint $\phi(z)\leq 0$ as the inclusion $\phi(z)\in \mathbb R^b_-$, we can define the linearized cone of (P) at $\bar z$ by $$\mathbb{L} (\bar z):=\{u\in\mathbb R^a|D\phi(\bar z)(u)\cap T_{\mathbb R^b_-}(\phi(\bar z))\neq\emptyset\}$$ and  the critical cone of (P) at $\bar z$ by} $$C(\bar z):=\{u\in \mathbb{L}(\bar z)| \nabla\varphi(z)u\leq0\}=\{u\in\mathbb R^a|D\phi(\bar z)(u)\cap T_{\mathbb R^b_-}(\phi(\bar z))\neq\emptyset, \nabla\varphi(\bar z)u\leq0\}.$$
By the directional differentiability of $\phi_j(1<j\leq b)$, one can easily obtain  $$D\phi(\bar z)(u)=D\phi_1(\bar z)(u)\times\Pi_{j=2}^b\phi'_j(\bar z;u),$$ and by the directional Lipschitz continuity of $\phi_1$ at $\bar z$ in direction $u$, $$D\phi_1(\bar z)(u)\cap\mathbb R_-\neq\emptyset\Leftrightarrow{(\phi_1)'_-(\bar z;u)\leq0}.$$ Hence the following equivalent definitions of $\mathbb{L}(\bar z)$ and $C(\bar z)$ hold:
\begin{align*}
	&\mathbb{L}(\bar z)=\{u\in\mathbb R^a|(\phi_j)_-'(\bar z;u)\leq0, j\in I_\phi(\bar z)\},\\
	&C(\bar z)=\{u\in L(\bar z)| \nabla\varphi(\bar z)u\leq0\}=\{u\in\mathbb R^a|(\phi_j)'_-(\bar z;u)\leq0, j\in I_\phi(\bar z), \nabla\varphi(\bar z)u\leq0\}.
\end{align*}

The concept of (Clarke) calmness for a mathematical program is first defined by Clarke \cite[Definition 6.41]{Clarke} and extended to a directional version by   Bai and Ye \cite{BY}. 
\begin{defn}[Directional Clarke Calmness]\cite[Definition 3.1]{BY}
	Suppose $\bar z$ is a local solution of {\rm (P)}. We say that {\rm (P)} is (Clarke) calm at $\bar z$ in direction $u$ if there exist positive scalars $\epsilon,\delta,\rho$, such that for any $\alpha \in\epsilon\mathbb B$ and any $z\in\bar z+{\cal V}_{\epsilon,\delta}(u)$ satisfying $\phi(z)+\alpha\leq 0$ one has,
	\begin{equation*}
		\varphi(z)-\varphi(\bar z)+\rho\|\alpha\|\geq0.
	\end{equation*}
\end{defn}
Similar to its nondirectional version, provided the objective function is directionally Lipschitz continuous, the directional calmness of problem (P) holds under the directional MSCQ, hence can be verified by directional constraint qualifications listed in Definition \ref{qp}.
\begin{lemma}\cite[Lemma 3.1]{BY}
	Let $\bar z$ solve $({\rm P})$ and $\varphi(z)$ be Lipschitz continuous at $\bar z$ in direction $u$. Suppose  that the system  $\phi(z)\leq 0 $ satisfies MSCQ  at $\bar z$ in direction $u$. Then $({\rm P})$ is calm at $\bar z$ in direction $u$.
\end{lemma}
 
\begin{prop}
	Let $\phi(\bar z)\leq0$ and suppose that $ \phi$ is directionally Lipschitz at $\bar z$ in direction {$u\in \mathbb{L}(\bar z)$.}
If the directional quasi-normality holds at $\bar z$ in direction $u$ for the inequality system $\phi(z)\leq 0$. Then the system $\phi(z)\leq0$ satisfies the directional MSCQ at $\bar z$ in direction $u$.
\end{prop}
\beginproof
The proof can be completed by following a similar process of the proof of \cite[Proposition 3.1]{BY}.
\endproof

We now summerize our discussions by the following implications under the Lipschitz continuity of the objective function in direction $u$:
\begin{center}
	Quasi-normality in direction $u$\\
	$\Downarrow$\\
	MSCQ in direction $u$\\
	$\Downarrow$\\
	Calmness in direction $u$.	
\end{center}

The following theorem shows that under the directional calmness of (P), a directional KKT condition holds at a local minimizer. This theorem   generalizes  \cite[Theorem 3.1]{BY} to allow one of the constraint functions to be not directionally differentiable. 
\begin{thm}\label{dKKT}
	Let $\bar{z}$ be a local minimizer of {\rm (P)} and $u\in C(\bar z)$. Suppose that $\varphi(z)$ is continuously differentiable at $\bar z$,  $\phi(z)$ is directionally Lipschitz continuous and $\phi_i(z)$ for $i=2,\ldots,b$ are directionally differentiable at $\bar z$ in direction $u$. 
	Suppose that Problem (P) is calm at $\bar z$ in direction $u$. Then there exist a multiplier $\lambda\in\mathbb R^{b}$ and  $\xi\in D\phi_1(\bar z)(u)$ such that
	$ 0\leq\lambda\perp \phi(\bar z),\ \lambda\perp d:=(\xi,\phi_2'(\bar z;u),\ldots,\phi_b'(\bar z;u))$ and 
	\begin{eqnarray*}
	0\in\nabla \varphi(\bar z)+\partial\langle\lambda,\phi\rangle(\bar z;u).
	\end{eqnarray*} 
\end{thm}
\beginproof
Since (P) is calm at $\bar z$ in direction $u$, there exist positive scalars $\epsilon,\delta,\rho$ such that
\begin{equation}\label{penalizedP}
	\varphi(z)+\rho\|\phi_+(z)\|\geq \varphi(\bar z)\qquad  \forall z\in \bar z+{\cal V}_{2\epsilon,\delta}(u).
\end{equation} 
Since $u\in C(\bar z)$, there exist $\xi\in D\phi_1(\bar z)(u)$ and sequences $t_k\downarrow0, u^k\rightarrow u$ such that $\lim_k(\phi_1(\bar z+t_ku^k)-\phi_1(\bar z))/t_k=\xi\leq0$. Then for sufficiently large $k$ we have  $$\phi(\bar z)+t_k(\xi,\phi_2'(\bar z;u),\ldots,\phi_b'(\bar z;u))\leq0$$ and hence
$$0\leq \frac{\|\phi_+(\bar z+t_ku^k)\|}{t_k}\leq \frac{\|\phi(\bar z+t_ku^k)-\phi(\bar z)-t_k(\xi,\phi_2'(\bar z;u),\ldots,\phi_b'(\bar z;u))\|}{t_k} .$$
It follows that $\lim_{k\downarrow \infty} (\|\phi_+(\bar z+t_ku^k)\|)/t_k=0$.

Since  $\bar z+t_k u^k\in \bar z+\operatorname{cl}({\cal V}_{\epsilon,\delta}(u))$ for $k$ sufficiently large,   by (\ref{penalizedP}),  $$\varphi(\bar z+t_ku^k)+\rho\|\phi_+(\bar z+t_ku^k)\|\geq\varphi(\bar z).$$ Together with 
$\nabla \varphi(\bar z)u\leq0$ we have
\begin{equation}\label{speed}
	\lim_{k\downarrow\infty}\frac{\varphi(\bar z+t_ku^k)+\rho\|\phi_+(\bar z+t_ku^k)\|-\varphi(\bar z)}{t_k}=0.
\end{equation}
The rest of the proof then follows from that of \cite[Theorem 3.1]{BY} by replacing $\frac{u}{k}$ and $\phi'(\bar z;u)$ with $t_ku^k$ and $(\xi,\phi_2'(\bar z;u),\ldots,\phi_b'(\bar z;u))$, respectively. For completeness we sketch the proof.
For each $k=0,1,\ldots$, define $\sigma_k:=2(\varphi(\bar z+t_ku^k)+\rho\|\phi_+(\bar z+t_ku^k)\|-\varphi(\bar z))$. 
If $\sigma_k\equiv0$, then for each large enough $k$, by (\ref{penalizedP}), $\bar z+t_ku^k$ is a global minimizer of the function $\varphi(z)+\rho\|\phi_+(z)\|+\delta_{\bar z+\operatorname{cl}({\cal V}_{\epsilon,\delta}(u))}(z)$. Since for each large enough $k$, $\bar z+t_ku^k$ is an interior point of $\bar z+\operatorname{cl}({\cal V}_{\epsilon,\delta}(u))$, by the well-known Fermat's rule and the calculus rule ({see e.g., \cite[Corollary 10.9]{RW}}), 
\begin{equation}\label{opk1}
	0\in\nabla\varphi(\bar z+t_ku^k)+\rho\partial(\|\phi_+{(\cdot)}\|)(\bar z+t_ku^k).
\end{equation}
Otherwise, without loss of generality, we assume that  for all $k$, $\sigma_k>0$.  Then by definition of $\sigma_k$,  the equality (\ref{speed}),
we have for $k$ sufficiently large, $$\varphi(\bar z+t_ku^k)+\rho\|\phi_+(\bar z+t_ku^k)\|+\delta_{\bar z+\operatorname{cl}({\cal V}_{\epsilon,\delta}(u))}(\bar z+t_ku^k)<\varphi(\bar z)+\sigma_k.$$   Define $\lambda_k:=\frac{2\|u^k\|rt_k}{\epsilon}\sqrt{\frac{\sigma_k\epsilon}{2\|u^k\|rt_k}}$ { with $r>0$}. 
Then by Ekeland's variation principle (see e.g., \cite[Theorem 2.26]{Aub2}),  there exists 
	$$\tilde z^k\in\operatorname{argmin}_{z}\left\{ \varphi(z)+\rho\|\phi_+( z)\|+\delta_{\bar z+\operatorname{cl}({\cal V}_{\epsilon,\delta}(u))}(z)+\frac{\sigma_k}{\lambda_k}\|z-(\bar z+t_ku^k)\|\right\}$$  satisfying 
	that $\|\tilde z^k-(\bar z+t_ku^k)\|\leq\lambda_k$. By (\ref{speed}), we can show that  for sufficiently small $r$,
$\tilde z^k$ is in the interior of $\bar z+\operatorname{cl}({\cal V}_{\epsilon,\delta}(u))$. Then by the well-known Fermat's rule, we obtain
\begin{equation}\label{opk2}
	0\in\nabla\varphi(\tilde z^k)+\rho\partial(\|\phi_+(\cdot)\|)(\tilde z^k)+\frac{\sigma_k}{\lambda_k}\bar{\mathbb B}.
\end{equation}
Using the  directionally Lipschitz continuity of  $\phi$  at $\bar z$ in direction $u$, by applying the chain rule for limiting subdifferential \cite[Corollary 3.43]{Aub2}, we have  $$\partial(\|\phi_+(\cdot)\|)(\tilde z^k)\subseteq\bigcup_{\zeta'\in  \partial\|(\cdot)_+\|(\phi(\tilde z^k))}\partial\langle\zeta',\phi\rangle(\tilde z^k).$$
Therefore by  $(\ref{opk1})$ or $(\ref{opk2})$, $\exists\zeta^k\in\partial\| (\cdot)_+\|(\phi(\bar z+t_ku^k))$ or $\exists\zeta^k\in\partial\| (\cdot)_+\|(\phi(\tilde z^k))$ such that
\begin{equation}\label{opk}
	0\in\nabla\varphi(\bar z+t_ku^k)+\rho\partial\langle\zeta^k,\phi\rangle(\bar z+t_ku^k),\ \mbox{or}\
	0\in\nabla\varphi(\tilde z^k)+\rho\partial\langle\zeta^k,\phi\rangle(\tilde z^k)+\frac{\sigma_k}{\lambda_k}\bar{\mathbb B}.
\end{equation}
By the Lipschitz continuity of the function $\|x_+\|$,
$\{\zeta^k\}$ is bounded. Without loss of generality, assume $\zeta:=\lim_k\zeta^k$. We also have $\lim_k(\bar z+t_ku^k-\bar z)/t_k=\lim_k(\tilde z^k-\bar z)/t_k=u$. Since $\sigma_k=o(t_k),\ \lim_k\frac{\sigma_k}{\lambda_k}=0$. Taking the limit of $(\ref{opk})$ as $k\rightarrow\infty$,   by \cite[Theorem 5.4]{Long} we have
$$0\in \nabla \varphi(\bar z)+\rho \partial \langle \zeta, \phi \rangle (\bar z; u).$$
Moreover by  \cite[Corollary 4.2]{BHA}, $$\zeta\in\partial(\|(\cdot)_+\|)(\phi(\bar z);(\xi,\phi_2'(\bar z;u),\ldots,\phi_b'(\bar z;u)))\subseteq N_{\mathbb R^b_-}(\phi(\bar z);(\xi,\phi_2'(\bar z;u),\ldots,\phi_b'(\bar z;u))).$$
The desired result holds by taking $d:=(\xi,\phi_2'(\bar z;u),\ldots,\phi_b'(\bar z;u))$ and $$\lambda_{\phi}:=\rho\zeta\in N_{\mathbb R^b_-}(\phi(\bar z);d)=\{\eta\in\mathbb R^b|0\leq\eta\perp\phi(\bar z), \eta\perp d\}.$$ 
\endproof
\section{Directional sensitivity analysis of the  Moreau envelope function}

Since the  Moreau envelope function $v_\gamma(x,y)$ is the optimal value function for  the following parametric mathematical program with parameter $(x,y)$:
\begin{align*}
\mbox{(P$^\gamma_{x,y}$)}~~~~\quad \min_w~ &\nu(x,y,w):=f(x,w)+\frac{1}{2\gamma}\|w-y\|^2\\
s.t.~ &g(x,w)\leq0,
\end{align*}
 in this section we study the sensitivity analysis of the  Moreau envelope function from the perspective that it is the optimal value function for the above parametric program. 
Note that  the feasible map of the problem {\rm (P$^\gamma_{x,y}$)} is the same as the one for problem {\rm (P$_x$)}:
$$\widetilde {\mathcal F}(x,y):={\mathcal F}(x)=\{w\in\mathbb R^m|g(x,w)\leq0\}.$$ We define 
the optimal solution map/the proximal map by
	\begin{equation*}
	\begin{aligned}
		S_{\gamma}(x,y) & :=\displaystyle  \operatorname{argmin}_{w} \left\{\left. f(x,w) + \frac{1}{2\gamma}\|w - y\|^2 \right| g(x,w) \le 0  \right\}.
	\end{aligned}
\end{equation*}
In this section, unless otherwise specified, we assume that for fixed $x$, $f$ is $\rho_f$-weakly convex in $w$ and $g$ is quasi-convex in $w$. Then for any postive number $\gamma<\frac{1}{2\rho_f}$, for each fixed $(x,y)$, the objective function $\nu(x,y,w)$ of the minimization problem (P$^\gamma_{x,y}$) is strongly convex and so there always exists a unique solution for each $(x,y)$. In this section we assume  $ \gamma\in(0,1/2\rho_f)$. This implies that the proximal map $S_\gamma$ is always single-valued. Besides, the objective function $\nu(x,y,w)$ is coercive with respect to $w$, and hence it is level-bounded. It follows that the restricted inf-compactness (see e.g., \cite[Definition 3.8]{GLYZ}) always holds.

We define the Lagrange function of {\rm (P$^\gamma_{x,y}$)}  by
\[
{\cal L}(x,y,w;\lambda):=f(x,w)+\frac{1}{2\gamma}\|w-y\|^2+g(x,w)^T\lambda,
\]
and the set of Lagrange multipliers of {\rm (P$^\gamma_{x,y}$)} associated with $w\in{\cal F}(x)$ by
\[
\Lambda(x,y,w):=\{\lambda\in\mathbb R^p|\nabla_w{\cal L}(x,y,w;\lambda)=0,\ g(x,w)^T\lambda=0,\ \lambda\geq0\}.
\] It is easy to see that for any $y\in {\cal F}(x)$, the set of Lagrange multipliers of {\rm (P$^\gamma_{x,y}$)} coincides with the set of Lagrange multipliers of 
{\rm (P$_x$)}, i.e.,
$\Lambda(x,y,y)=\Lambda(x,y).$

 Let $(u,v)$ satisfy $\nabla g( x, w)(u,v) \in T_{\mathbb R^p_-}(g( x, w))=\{(\xi,\zeta) \mid \nabla g_i(x,w)(\xi,\zeta) \leq 0 \ \ i\in I_g(x,w)\}$.
We denote the set of directional Lagrange multipliers of {\rm (P$^\gamma_{x,y}$)} associated with $w\in{\cal F}(x)$ in direction $(u,v)$ by
\[
\Lambda(x,y,w;u,v):=\{\lambda\in\mathbb R^p_+|\nabla_w{\cal L}(x,y,w;\lambda)=0,\ g(x,w)^T\lambda=0,\ \lambda^T\nabla g(x,w)(u,v)=0\}.
\]

We divide our discussion into three cases: first, we show that sensitivity analysis of value functions in convex analysis can be extended to the case when $v_\gamma$ is weakly convex. 
In the second and third subsections, without weak convexity of $v_\gamma$, we obtain sensitivity analysis under FOSCMS and RCR regularity+RS, respectively.
\subsection{Sensitivity analysis under weak convexity}
Recently, Gao et al. \cite{GLYZ} studied weak convexity and sensitivity analysis of $v_\gamma$. Particularly, in \cite[Theorems 1-3]{GLYZ}, they show that $v_\gamma$ is weakly convex when $f$ is weakly convex and $g$ is convex. And from their proof, one can easily find that convexity of $g$ can be weakened to quasiconvexity. Moreover, in this paper we only need the local Lipschitz continuity of the function $v_\gamma(x,y)$.
\begin{prop}\label{wc}
	Assume that $g(x,y)$ is quasiconvex on  $ \mathbb R^n\times\mathbb R^m$, $f$ is $\rho_f$-weakly convex  with $\rho_f>0$ on $\operatorname{gph}\mathcal F$. Then the Moreau envelope function $v_\gamma(x,y)$ is $\rho_v$-weakly convex for any $\gamma\in(0,1/\rho_f)$ satisfying $\rho_v\geq\frac{\rho_f}{1-\gamma\rho_f}$. If $\mathcal F(x)\neq\emptyset$ for all $x$ in a neighborhood of $\bar x\in\mathbb R^n$, then $v_\gamma(x,y)$ is Lipschitz continuous around $(\bar x,y)$ for any $y\in\mathbb R^m$.
\end{prop}

Under the assumptions in Proposition \ref{wc}, Gao et al. \cite{GLYZ} presented lower estimates for the subdifferential of $v_\gamma$. We now extend their results to directional subdifferential defined in Definition \ref{ads}.
\begin{prop}\label{wcsbdiff}
	Assume that $g(x,y)$ is quasiconvex on $\mathbb R^n\times\mathbb R^m$,  $f$ is $\rho_f$-weakly convex  with $\rho_f>0$ on $\operatorname{gph}\mathcal F$. 
	Let $\bar w\in S_\gamma(\bar x,\bar y)$ and suppose that  $\mathcal F(x)\neq\emptyset$ for all $x$ in a neighborhood of $\bar x$. Then $v_\gamma(x,y)$ is Lipschitz continuous  and  directionally differentiable at $(\bar x, \bar y)$.
	One has
	\begin{equation}\label{wclowerest1}
			\partial v_\gamma(\bar x,\bar y)\supseteq
		\bigcup_{\lambda \in\Lambda(\bar x,\bar y,\bar w)}
		\left(	\left\{
		\nabla_xf(\bar x,\bar w)+\nabla_xg(\bar x,\bar w)^T\lambda
		\right\}\times\{(\bar y-\bar w)/\gamma\}\right),
	\end{equation}
and the equality holds if Guignard constraint qualification  (Guignard CQ) for the system $g(x,w)\leq 0$  holds at $(\bar x,\bar w)$, i.e.,
\begin{equation}\label{Guignard}
	N_{\operatorname{gph}\mathcal F}(\bar x,\bar w)=\left\{
	\sum_{i=1}^l\lambda_i\nabla g_i(\bar x,\bar w)|\lambda\geq0,\sum_{i=1}^l\lambda_ig_i(\bar x,\bar w)=0
	\right\}.
\end{equation} If Guignard CQ (\ref{Guignard}) holds, then for any direction $(u,v)\in\mathbb R^n\times\mathbb R^m$ satisfying $\nabla g(\bar x,\bar w)(u,v)\in T_{\mathbb R^p_-}(g(\bar x,\bar w))$, we have
	\begin{equation}\label{wclowerest2}
	\partial v_\gamma(\bar x,\bar y;u,v)=
	\bigcup_{\lambda\in{\Lambda(\bar x,\bar y,\bar w;u,v)}}
	\left(	\left\{
	\nabla_xf(\bar x,\bar w)+\nabla_xg(\bar x,\bar w)^T\lambda
	\right\}\times\{(\bar y-\bar w)/\gamma\}\right).
\end{equation}
	
	\begin{equation}\label{wcdird}
		v'_\gamma(\bar x,\bar y;u,v)=\max_{\lambda\in{\Lambda(\bar x,\bar y,\bar w;u,v)}}\nabla\mathcal L(\bar x,\bar y,\bar w;\lambda)(u,v).
	\end{equation}
\end{prop}
\beginproof
Since $\mathcal F(x)\neq\emptyset$ for all $x$ in a neighborhood of  $\bar x$, by Proposition \ref{wc}, $v_\gamma$ is weakly convex and 
Lipschitz continuous around $(\bar x,\bar y)$. 
The first estimate (\ref{wclowerest1}) follows by \cite[Theorem 5]{GLYZ}. We now prove formula (\ref{wclowerest2}). By Proposition \ref{Prop2.2}, $(\xi,\zeta)\in\partial v_\gamma(\bar x,\bar y;u,v)$ if $(\xi,\zeta)\in\partial v_\gamma(\bar x,\bar y)$ and $\langle (\xi,\zeta),(u,v)\rangle=\sigma_{\partial v_\gamma(\bar x,\bar y)}(u,v)$, 
where  $\sigma_C$ denotes the support function of set $C$. When the equality in (\ref{wclowerest1}) holds, we have
\begin{align*}
\sigma_{\partial v_\gamma(\bar x,\bar y)}(u,v)&=\max_{\lambda\in\Lambda(\bar x,\bar y,\bar w)}
\left(	
\nabla_xf(\bar x,\bar w)u+\sum^l_{i=1}\lambda_i\nabla_xg_i(\bar x,\bar w)u
+ v^T(\bar y-\bar w)/\gamma \right).
\end{align*} 
Now suppose $\nabla g(\bar x,\bar w)(u,v)\in T_{\mathbb R^p_-}(g(\bar x,\bar w))$ and Guignard CQ (\ref{Guignard}) holds. Then  for any vector $\lambda \in\Lambda(\bar x,\bar y,\bar w)$, we have $\lambda \in N_{\mathbb R^p_-}(g(\bar x,\bar w))$. Therefore $\lambda^T\nabla g(\bar x,\bar w)(u,v)\leq0$. This means $\langle (\xi,\zeta),(u,v)\rangle=\sigma_{\partial v_\gamma(\bar x,\bar y)}(u,v)$ if and only if $\xi=\nabla_xf(\bar x,\bar w)+\sum^l_{i=1}\lambda_i\nabla_xg_i(\bar x,\bar w), \lambda\in\Lambda(\bar x,\bar y,\bar w)$ and $\lambda^T\nabla_xg(\bar x,\bar w)(u,v)=0$, that is, $\lambda\in\Lambda(\bar x,\bar y,\bar w;u,v)$. Hence, when Guignard CQ holds, one has the formula (\ref{wclowerest2}). 
(\ref{wcdird}) follows by applying Proposition \ref{Prop2.2}.
\endproof

{
It is known that Guignard CQ holds under  MSCQ. The reader can refer to \cite[Proposition 5]{YYZZ} for more sufficient conditions for Guignard CQ.
}

\subsection{Sensitivity analysis without weak convexity}
First we state a result on the directional Lipschitz continuity of $v_\gamma$.
{\begin{prop}\cite[Theorem 4.1]{BY}\label{Lip}
	Let $\bar w\in S_\gamma(\bar x,\bar y)$. Suppose that the feasible map ${\cal F}$ satisfies RS at  $(\bar x, \bar w)$ in direction $u$.
	Then $v_\gamma(x,y)$ is Lipschitz continuous at $(\bar x,\bar y)$ in direction $(u,v)$ for any $v$. 
\end{prop}
\beginproof
 By Definition \ref{dRS}, ${\cal F}(x)$ satisfying RS at  $(\bar x, \bar w)$ in direction $u$ is equivalent to $\widetilde{\cal F}(x,y)$ satisfying RS at  $(\bar x, \bar y, \bar w)$ in direction $(u,v)$ for any $v$. By \cite[Theorem 4.1]{BY} $v_\gamma(x,y)$ is Lipschitz continuous at $(\bar x,\bar y)$ in direction $(u,v)$ for any $v$. 
\endproof}

 Observe that $\nu(x,y,w)=f(x,w)+\frac{1}{2\gamma}\|w-y\|^2$ is smooth with the gradient equal to $$\nabla\nu(x,y,w)=(\nabla_xf(x,w),(y-w)/\gamma,\nabla_wf(x,w)+(w-y)/\gamma).$$ So by the directional Danskin's theorem \cite[Proposition 4.4]{BY2023}, we conclude that the Clarke subdifferential of $v_\gamma( x,y) $ with respect to variable $y$  is a singleton. Hence
$v_\gamma(x,y)$ is continuously differentiable in variable $y$ and $\nabla_y v_\gamma (\bar x,\bar y)=(\bar y-\bar w)/\gamma$ with $\bar w\in S_\gamma(\bar x,\bar y)$. 
 We have
\begin{align}\label{ddformula}
	(v_\gamma)_+' (\bar x,\bar y;u,v)&=(v_\gamma)_+'(\bar x,\bar y;u,0)+v^T(\bar y-\bar w)/\gamma,\notag\\
	(v_\gamma)_-' (\bar x,\bar y;u,v)&=(v_\gamma)_-'(\bar x,\bar y;u,0)+v^T(\bar y-\bar w)/\gamma,\notag\\
	v_\gamma' (\bar x,\bar y;u,v)&=v_\gamma'(\bar x,\bar y;u,0)+v^T(\bar y-\bar w)/\gamma,
\end{align}
where $\bar w\in S_\gamma(\bar x,\bar y)$.

To obtain directional differentiability and subdifferentiability of $v_\gamma$ with respect to variable $x$, we need to assume some regularity conditions. In the following subsections, we study the sensitivity analysis of $v_\gamma$ under  FOSCMS,  and RCR regularity$+$RS, respectively.

\subsubsection{Sensitivity analysis under FOSCMS}

Since the solution map $S_\gamma(x,y)$ is always uniformly compact around $(\bar x,\bar y)$, by \cite[Theorem 4.3 and Corollaries 4.7 and 4.8]{BY2023-1}, under MFCQ or FOSCMS we have the following estimates for Dini directional derivatives of $v_\gamma(x,y)$. 
\begin{prop}\label{dini}  Let $\bar w\in S_\gamma(\bar x,\bar y)$. Suppose either MFCQ holds for the system $g(\bar x, w)\leq0$ at $w=\bar w$ or  one of the following FOSCMS holds:
\begin{itemize} 
\item Suppose for any direction $d\in {\rm Im} D_xg(\bar x,\bar w; u)$, there exists $v\in \mathbb{R}^m$
satisfying 
\begin{equation}\label{eqnImnew}
d_i+\nabla_y g_i(\bar x,\bar w)v\leq0,\ i\in I_g(\bar x,\bar w).\end{equation}
Suppose that FOSCMS for the system $g(\bar x, w)\leq0$ at $w=\bar w$ holds in direction $d+\nabla_y g(\bar x,\bar w)v$, i.e.,
\begin{eqnarray*}
	\left\{
	\begin{array}{ll}
		\nabla_yg(\bar x,\bar w)^T\lambda=0,\\
		0\leq\lambda\perp g(\bar x,\bar w),\ \lambda\perp (d+\nabla_y g(\bar x,\bar w)v )
	\end{array}
	\right.\implies \lambda=0,
\end{eqnarray*}
 for each $(d,v)\not =0$ satisfying (\ref{eqnImnew}) and $d\in {\rm Im} D_xg(\bar x,\bar w; u)$. Moreover suppose $\mathbb L(\bar x,\bar w;\pm u)$ 
is nonempty and FOSCMS holds at $\bar w$ in direction $\nabla g(\bar x,\bar w)(\pm u,v)$ for any $v\in\mathbb L(\bar x,\bar w;\pm u)$. 
 \item $\nabla_x g(\bar x,\bar y)u\not =0$, $ \mathbb L(\bar x,\bar w;\alpha u)$ is nonempty for any $\alpha \in \mathbb{R}$ and FOSCMS holds at $\bar w$ in direction $\nabla g(\bar x,\bar w)(\alpha u, v)$ for any $(\alpha, v)\not =(0,0)$ with $v\in \mathbb L(\bar x,\bar w;\alpha u)$. 
 \end{itemize} 
  Then 
	 for any direction $v\in \mathbb{R}^m$,
	 \begin{equation}\label{neweqn20}
	\begin{aligned}
		\min_{\lambda\in  \Lambda(\bar x, \bar y,\bar w)} \nabla_x \mathcal L(\bar x,\bar y,\bar w;\lambda)u+v^T(\bar y-\bar w)/\gamma\leq &(v_\gamma)_-'(\bar x,\bar y;u,v)  
		\leq (v_\gamma)_+'(\bar x,\bar y;u,v) \\ \leq 
		{\max_{\lambda\in  \Lambda(\bar x, \bar y,\bar w)}}\nabla_x \mathcal L(\bar x,\bar y,\bar w;\lambda)u+v^T(\bar y-\bar w)/\gamma.
	\end{aligned} 
	\end{equation}
	Furthermore if  $\Lambda(\bar x, \bar y, \bar{w})=\{\bar\lambda\}$ is a singleton, 
	then $v_\gamma(x,y)$ is Hadamard directionally differentiable and 
	$$(v_\gamma)'(\bar x,\bar y; u, v)= 
	\nabla_x \mathcal L(\bar x,\bar y,\bar w;\bar \lambda)u+ v^T(\bar y-\bar w)/\gamma.$$
\end{prop} 
\beginproof
When MFCQ holds, the proof is complete from \cite[Corollary 4.7]{BY2023-1}. Next, we prove the results under FOSCMS. As discussed at the beginning of Section 4.2,  under the assumption of this section, the restricted inf-compactness condition holds automatically. By Proposition \ref{Prop2.4new}, under the assumed conditions $\mathcal F$ satisfies RS in direction $(u,v)$. The estimate (\ref{neweqn20}) follows from \cite[Theorem 4.3]{BY2023-1}.
\endproof

We now study the subdifferentiability of $v_\gamma$. Let $w\in {\cal F}(x)$. 
		For direction $(u,d) \in \mathbb{R}^n \times \mathbb{R}^m$, define
		\begin{eqnarray*}
			W( x, y, w; u,d):= \bigcup_{\lambda\in \Lambda(x,y,w;u,d)} \left\{(\nabla_x f(x, w)+\nabla_x g(x,w)^T\lambda)\right\} \times \{(y-w)/\gamma\}.
		\end{eqnarray*}
		
		The following inner calmness* property will be used in this section to obtain simplified estimates for directional subdifferentials of $v_\gamma$.
		\begin{defn}[Directional Inner Calmness*](Benko et al. \cite{B2021}) The optimal solution map $S_\gamma(x,y)$ is said to be inner calm* at $(\bar x,\bar y,\bar w)$ in direction $(u,v)$ if there exists $\kappa>0$ such that for every sequence $(x^k,y^k)\xrightarrow{(u,v)}(\bar x,\bar y)$, 
			there exist a subsequence $K$ of the set of nonnegative integers $\mathbb{N}$, 
			{together with a sequence  $ w^k\in S_\gamma(x^k,y^k)$ for $k\in K$ such that} 
			\begin{equation*}
				\| w^k-\bar w\|\leq\kappa\| (x^k,y^k)-(\bar x,\bar y)\|.
			\end{equation*} 
			Particularly, we say $S_\gamma(x,y)$ is inner calm* at $(\bar x,\bar y,\bar w)$ if $(u,v)=(0,0)$.
		\end{defn}

		Based on \cite[Proposition 4.2]{BY2023} we can obtain the following results.   
		\begin{prop}\label{est1}  Let $\bar w\in S_\gamma(\bar x,\bar y)$ and $u\in \mathbb{R}^n$. 
Suppose that the system $g(x,w)\leq0$ satisfies MSCQ at  $(\bar x,\bar w)$ and $\mathcal F$ satisfies RS at $(\bar x,\bar w)$ in direction $u$.
Then $v_\gamma(x,y)$ is Lipschitz around $(\bar x,\bar y)$ in direction $(u,v)$, and
\begin{align}\label{upperestimate}
\emptyset \not = &\partial v_\gamma(\bar x,\bar y;u,v) \subseteq \bigcup_{d\in\mathcal C(\bar x,\bar y,\bar w;u)\cup(\mathcal C(\bar x,\bar y,\bar w;0)\cap\mathbb S)} W(\bar x,\bar y,\bar w;u,d),
\end{align}
where
\begin{align*}
&\mathcal C(x,y,w;u)\\
:=&\left \{d ~\Bigg| \begin{array}{l}
	(v_\gamma)_-'(x,y;u,0) \leq   \nabla f(x,w)(u,d)+ {d^T(w-y)/\gamma } \leq (v_\gamma)_+'(x,y;u,0)\\\nabla g_i(x,w)(u,d)\leq 0, \ \ \forall i\in I_g(x,w) \end{array} 
\right \}.
\end{align*} 
				Moreover if $S_\gamma(x,y)$ is inner calm* at $(\bar x,\bar y,\bar w)$ in direction $(u,v)$, then the upper estimate (\ref{upperestimate}) can be replaced by the smaller set
$\bigcup_{d \in\mathcal C(\bar x,\bar y,\bar w;u)} W(\bar x,\bar y,\bar w;u,d)$.				
%
		\end{prop}
	\beginproof  As discussed at the beginning of Section 4.2,  under the assumption of this section, the restricted inf-compactness condition holds automatically. As $S_\gamma(\bar x,\bar y)=\{\bar w\}$ and the system $g(x,w)\leq0$ satisfies MSCQ at $(\bar x,\bar w)$, by \cite[Theorem 3.1]{BY2023} the upper estimate (\ref{upperestimate}) holds. Besides, since $\mathcal F$ satisfies RS at $(\bar x,\bar w)$ in direction $u$, the directional Lipschitz continuity of $v_\gamma$ along direction $(u,v)$ follows directly from Proposition \ref{Lip}. It follows that $\partial v_\gamma(\bar x,\bar y;u,v) $ is nonempty.
\endproof

{Since by Proposition \ref{Prop2.4new}, FOSCMS in certain directions or MFCQ for the system $g(\bar x, w)\leq 0$ holding at $\bar w$ implies RS, we can replace RS by FOSCMS/MFCQ and obtain the following corollary from Propositions \ref{dini} and \ref{est1}.
\begin{cor}\label{cor4.1}  Let $\bar w\in S_\gamma(\bar x,\bar y)$ and   $u\in \mathbb{R}^n$. Then under one of the following conditions,  $v_\gamma(x,y)$ is directionally differentiable and Lipschitz around $(\bar x,\bar y)$ in direction $(u,v)$ for any $v\in \mathbb{R}^m$, and (\ref{neweqn20}), (\ref{upperestimate}) hold.
\begin{itemize}
\item[(a)] Suppose $\nabla_x g(\bar x,\bar w)u\not =0$ and $\mathbb L(\bar x,\bar w;\alpha u)\not =\emptyset$ for any $\alpha\in\mathbb R$. Suppose  that the system $g(x,w)\leq0$ satisfies MSCQ at  $(\bar x,\bar w)$ and FOSCMS for  the system $g(\bar x, w)\leq0$ at $w=\bar w$ holds in direction $\nabla g(\bar x,\bar w)(\alpha u, v)$ for any $(\alpha,v)\neq0$ with $v\in\mathbb L(\bar x,\bar w;\alpha u)$.
\item[(b)] Suppose  that the system $g(x,w)\leq0$ satisfies MFCQ at  $(\bar x,\bar w)$
\end{itemize}
\end{cor}
\beginproof
(a) When $\nabla_xg(\bar x,\bar w)u\neq0$, $\{\alpha\nabla_xg(\bar x,\bar w)u|\alpha\geq0\}= {\rm Im}D_xg(\bar x,\bar w;u)$. Then $d\in{\rm Im}D_xg(\bar x,\bar w;u)$ if and only if $d=\nabla g(\bar x,\bar w)(\alpha u)$ for some $\alpha\geq0$. Hence condition (\ref{eqnImnew}) reduces to the nonemptyness of $\mathbb L(\bar x,\bar w;\alpha u)$ for all $\alpha\geq0$. {And FOSCMS at $\bar w$ in direction $d+\nabla_y g(\bar x,\bar w)v$ for each $(d,v)\not =0$ satisfying (\ref{eqnImnew}) and $d\in {\rm Im} D_xg(\bar x,\bar w;u)$ is equivalent to FOSCMS at $\bar w$ in direction $\nabla g(\bar x,\bar w)(\alpha u,v)$ for any $(\alpha, v)\not =(0,0)$ satisfying $v\in\mathbb L(\bar x,\bar w;\alpha u)$}. Estimate (\ref{neweqn20}) follows from applying Proposition \ref{dini} immediately. And (\ref{upperestimate}) holds because the directional RS follows from Proposition \ref{Prop2.4new} and MSCQ for the system $g(x,w)\leq0$ is assumed. 

Now we prove (b). By Proposition \ref{dini}, estimate (\ref{neweqn20}) holds under MFCQ. Furthermore, MFCQ implies RS of $\mathcal F$ at $(\bar x,\bar w)$ by Proposition \ref{Prop2.4new}, hence implies MSCQ of system $g( x,w)\leq0$ and directional RS of $\mathcal F$. Then by Proposition \ref{est1}, the estimate (\ref{upperestimate}) holds.
\endproof}
\subsection{Sensitivity analysis under RCRCQ+RS}
It is well-known in nonlinear optimization that RCR-regularity is not comparable to MFCQ. Now we apply directional sensitivity analysis from \cite{BY} and obtain the following directional differentiability results of $v_\gamma(x,y)$. 

\begin{prop}\label{ddv}
	Let $\bar w\in S_\gamma(\bar x,\bar y)$. Assume that ${\cal F}$ is RCR-regular at  $(\bar x, \bar w)$ and RS is satisfied at  $(\bar x, \bar w)$ in direction $u$.  Then the function $v_\gamma(x,y)$ is directionally Lipschitz and directionally differentiable at $(\bar x, \bar y)$ in  direction $(u, v)$ for any $v$ and		
	\begin{align}
		&v_\gamma'(\bar x,\bar y;u,v)=\max_{\lambda\in\Lambda(\bar x,\bar y,\bar w)}\nabla_x{\cal L}(\bar x,\bar y,\bar w;\lambda)u+v^T(\bar y-\bar w)/\gamma, \label{equality}\\
		 	&\emptyset \not = \partial v_\gamma(\bar x,\bar y;u,v)\subseteq	\partial^c v_\gamma(\bar x,\bar y;u,v) \subseteq W(\bar x, \bar y,\bar w; u,d) \quad \forall d\in \mathcal C(\bar x,\bar y,\bar w;u), \label{secondincl} 	\end{align} 
	where 
	\begin{align*}
		&\mathcal C(\bar x, \bar y,\bar w;u) = \left \{d ~\Bigg| \begin{array}{l}
			\max_{\lambda\in\Lambda(\bar x,\bar y,\bar w)}\nabla_x{\cal L}(\bar x,\bar y,\bar w;\lambda)u	 = \nabla f(\bar x,\bar w)(u,d)+ {d^T(\bar w-\bar y)/\gamma} \\ \nabla g_i(\bar x,\bar w)(u,d)\leq 0, \ \ \forall i\in I_g(\bar x,\bar w) \end{array} 
		\right \}.
	\end{align*} 

\end{prop}
\beginproof As discussed in the beginning of this section, the restricted inf-compactness always holds for problem {\rm (P$^\gamma_{x,y}$)} at any point $(\bar x,\bar y)$ and in any direction $(u,0)$. Since ${\cal F}(x)$ is independent of variable $y$, ${\cal F}$ satisfying  RS  at  $(\bar x, \bar w)$ in direction $u$ is equivalent to saying that ${\cal F}$ satisfies RS  at  $(\bar x, \bar y, \bar w)$ in direction $(u,0)$. Hence by \cite[Corollary 4.1]{BY2023-1}, $v_\gamma(x,y)$ is continuous in direction $(u,0)$ and  $S_\gamma(\bar x,\bar y;u,0)\not =\emptyset$, where $S_\gamma(\bar x,\bar y;u,0)$ denotes the directional solution in direction $(u,0)$ (see \cite[Definition 4.5]{BY}).  Hence $S_\gamma(\bar x,\bar y;u,0)=\{\bar w\}$.
By \cite[Proposition 4.3]{BY}, 
\[ 
v_\gamma'(\bar x,\bar y;u,0)=\max_{\lambda\in\Lambda(\bar x,\bar y,\bar w)}\nabla_x{\cal L}(\bar x,\bar y,\bar w;\lambda)u.
\]
Noting that $$v_\gamma'(\bar x,\bar y;u,v)=v_\gamma'(\bar x,\bar y;u,0)+\nabla_yv_\gamma(\bar x,\bar y)v\mbox{ and } \nabla_y v_\gamma (\bar x,\bar y)=(\bar y-\bar w)/\gamma,$$ the proof of (\ref{equality}) is completed. By \cite[Theorem 4.3]{BY}, since  $S_\gamma(\bar x,\bar y;u,v)=\{\bar w\}$ is a singleton, we have $ \partial^cv_\gamma(\bar x,\bar y;u,v) \subseteq \operatorname{co} W(\bar x, \bar y,\bar w; u,d) \quad \forall d \in d\in \mathcal C(\bar x,\bar y,\bar w;u)$. { Since $W(\bar x, \bar y,\bar w; u,d)$ is convex,} 
the second inclusion in (\ref{secondincl}) follows.
\endproof

\section{Necessary optimality conditions for problem {(VP)$_\gamma$}}

The main purpose of this section is to apply Theorem \ref{dKKT} to problem (VP)$_\gamma$ and the result of the directional sensitivity analysis of the Moreau envelope function in Section 4 to derive a verifiable necessary optimality condition for (VP)$_\gamma$ under a weak and verifiable constraint qualification. Recall that if the lower-level program is generalized convex, then the problem 
(VP) and the poblem (VP)$_\gamma$ are equivalent for any $ \gamma>0$. Thus under the lower-level generalized convexity assumption, we will provide a result for the original bilevel problem.  Besides, even without the generalized convexity,  an optimality condition of (VP)$_\gamma$ can still serve as an optimality condition for the original bilevel program. For example,  if the graph of the optimal solution map $S(x)$ and the one for the lower-level stationary map 
$$ \left  \{y|0\in \nabla_y f(x,y)+N_{{\cal F}(x)}(y)\right \}$$
 have the same localization, 
 then there is an equivalence of local optimality between (VP) and (VP)$_\gamma$ over a neighborhood of $(\bar x, \bar y)$. Consequently, necessary optimality conditions of (VP)$_\gamma$ can still be used for (VP).
Throughout this section, unless otherwise specified, we assume that for fixed $x$, $f$ is weakly convex in $y$ and $g$ is quasi-convex in $y$.



Since for any $(\bar x,\bar y)$
 we have
\[
(f-v_\gamma)'_-(\bar x,\bar y;u,v)=\nabla f(\bar x,\bar y)(u,v)-(v_\gamma)'_+(\bar x,\bar y;u,v)=\nabla f(\bar x,\bar y)(u,v)-(v_\gamma)'_+(\bar x,\bar y;u,v),
\] we can define the linearized cone of (VP)$_\gamma$ at a feasible solution $(\bar x,\bar y)$ by  
$$\mathbb{L}(\bar x,\bar y):=
\left \{(u,v) \Bigg| \begin{array}{l}
\nabla f(\bar x,\bar y)(u,v)\leq  (v_\gamma)_+'(\bar x, \bar y;u,v), \\
\nabla g_i(\bar x,\bar y)(u,v)\leq 0,\ \forall i\in I_{g}(\bar x,\bar y), {\nabla G_i}(\bar x,\bar y)(u,v)\leq 0,\ \forall i\in I_{G}(\bar x,\bar y)\end{array} \right \}.
$$


Let $(\bar x,\bar y)$ be a feasible solution of (VP)$_\gamma$. Then for all $y\in {\cal F}(x)$,
$$0=f(\bar x,\bar y)-v_\gamma (\bar x,\bar y) \leq f(x,y)-v_\gamma(x,y).  $$ Hence under the Abadie constraint qualification,  one always has $\nabla f(\bar x,\bar y)(u,v)\geq (v_\gamma)'_-(\bar x,\bar y;u,v) $ for all $(u,v)$ satisfying $\nabla g_i(\bar x,\bar y)(u,v)\leq 0,\ i\in I_{g}(\bar x,\bar y)$.      
The Abadie constraint qualification is very weak and {it would be satisfied under} the constraint qualifications we impose. We may denote  the critical cone  of (VP)$_\gamma$ at $(\bar x,\bar y) $ by
\begin{align}\label{def_C}
C(\bar x,\bar y)
:= & \{(u,v)\in \mathbb{L}(\bar x,\bar y) ~|~ \nabla F(\bar x,\bar y)(u,v)\leq0\} \notag\\
 = &	\left \{(u,v) ~\Bigg| ~\begin{array}{l}
	{(v_\gamma)'_-(\bar x, \bar y;u,v)} \leq \nabla f(\bar x,\bar y)(u,v)\leq  (v_\gamma)_+'(\bar x, \bar y;u,v),\\ \nabla F(\bar x,\bar y)(u,v)\leq0,\\
	\nabla g_i(\bar x,\bar y)(u,v)\leq 0,\ \forall i\in I_{g}(\bar x,\bar y), {\nabla G_i}(\bar x,\bar y)(u,v)\leq 0,\ \forall i\in I_{G}(\bar x,\bar y)\end{array} \right \}.
\end{align}
For $(\bar x,\bar y)\in \operatorname{gph}{\cal F}$ and $u\in \mathbb{R}^n$, define
\begin{align*}
					&\mathcal C(\bar x,\bar y;u)=\left \{v ~\Bigg| \begin{array}{l}
						(v_\gamma)'_-(\bar x,\bar y;u,v) \leq \nabla f(\bar x,\bar y)(u,v)\leq  (v_\gamma)_+'(\bar x, \bar y;u,v)\\ \nabla g_i(\bar x,\bar y)(u,v)\leq 0,\ \forall i\in I_{g}(\bar x,\bar y) \end{array} 
					\right \}.
				\end{align*} 
Define the Lagrange function of $(\mathrm{P_x})$  associated with $y\in{\cal F}(x)$ by 
\[
{\cal L}(x,y;\lambda):= {\cal L}(x,y,y;\lambda)=f(x,y)+g(x,y)^T\lambda,
\]
the set of Lagrange multipliers of $(\mathrm{P_x})$ by
\[
\Lambda(x,y):=\Lambda(x,y,y)=\{\lambda\in\mathbb R^p|\nabla_y{ \cal L}(x,y;\lambda)=0,\ g(x,y)^T\lambda=0,\ \lambda\geq0\}
\]
{and the set of directional multipliers associated with direction $(u,v)$ with $\nabla g( x, y)(u,v) \in T_{\mathbb R^p_-}(g( x, y))$ by
\[
\Lambda(x,y;u,v):=\{\lambda\in\mathbb R^p_+|\nabla_y{\cal L}(x,y;\lambda)=0,\ g(x,y)^T\lambda=0,\ \lambda^T\nabla g(x,y)(u,v)=0\}.
\]}



{It is known from \cite{YZ95}} that classical constraint qualifications such as the NNAMCQ fails at each feasible point of (VP). Recently, it is proved in \cite{BY} that the FOSCMS also fails at any feasible point of (VP) along any critical direction. Similar to \cite[Proposition 5.1]{BY} we can show that even the weaker constraint qualification FOSCMS fails at any feasible point of (VP)$_\gamma$ along any critical direction. As discussed in Section \ref{Sec3}, one can obtain a necessary optimality condition under the calmness condition.




In \cite[Definition 3.1]{YZ95}, the partial calmness condition was introduced for (VP). We now define a directional version of the partial calmness for our problem  {\rm (VP)$_\gamma$}.
\begin{defn}[Directional partial calmness]
	Suppose $(\bar x,\bar y)$ is a local solution of {\rm (VP)$_\gamma$}. We say that {\rm (VP)$_\gamma$} is partially calm at $(\bar x,\bar y)$ in direction $(u,v)$ if there exist positive scalars $\epsilon,\delta,\rho$, such that for any $(\alpha,\beta) \in\epsilon\mathbb B$ and any $(x,y)\in\operatorname{gph}\mathcal F\cap((\bar x,\bar y)+{\cal V}_{\epsilon,\delta}(u,v))$ satisfying $f(x,y)-v_\gamma(x,y)+\alpha\leq 0, G(x,y)+\beta\leq0$ one has,
	\begin{equation*}
		F(x,y)-F(\bar x,\bar y)+\rho\|(\alpha,\beta)\|\geq0.
	\end{equation*}
\end{defn}
It is easy to see that the partial calmness condition above is equivalent to exact penalization of the constraints $f(x,y)-v_\gamma(x,y)\leq 0 $ and $G(x,y)\leq 0$ in a directional neighborhood of $(u,v)$. It is easy to see that the combination of the directional partial calmness and directional MSCQ of the system $g(x,y)\leq0$ in direction $(u,v)$ can imply the directional calmness condition of (VP)$_\gamma$ in direction $(u,v)$.  Now we  
 apply Theorem \ref{dKKT} to (VP)$_\gamma$ to { obtain its  necessary optimalitys   condition under} the partial calmness condition. Note that under the lower-level generalized convexity assumption, the following theorem  provides a result for the original bilevel problem.

\begin{thm}\label{opt0}
	Let $(\bar x,\bar y)$ be a local minimizer of {\rm (VP)$_\gamma$}.
Suppose the value function  $v
	_\gamma$ is directionally Lipschitz continuous at $(\bar x,\bar y)$ in direction   $(u,v)\in C(\bar x,\bar y)$, defined in \eqref{def_C}. Moreover suppose {\rm (VP)}$_\gamma$ is either calm or partially calm  at $(\bar x,\bar y)$ in direction $(u,v)$ and MSCQ for the system $g(x,y)\leq 0$ at $(\bar x,\bar y)$ in direction $(u,v)$.
	Then there exists  $(\alpha, \lambda_g, \lambda_G)\in \mathbb R^{1+p+q}$ such that
	\begin{equation*}
	\begin{aligned}
		 0\in \nabla F(\bar x,\bar y)+\alpha \nabla f(\bar x,\bar y)-\alpha \partial^c v_\gamma(\bar x,\bar y;u,v) +
		\nabla g(\bar x,\bar y)^T\lambda_g+\nabla G(\bar x,\bar y)^T\lambda_G, \\
		\alpha\geq0,\ 0\leq\lambda_g\perp g(\bar x,\bar y),\ \lambda_g\perp\nabla g(\bar x,\bar y)(u,v),\ 0\leq\lambda_G\perp G(\bar x,\bar y),\ \lambda_G\perp\nabla G(\bar x,\bar y)(u,v).
	\end{aligned}
	\end{equation*}
\end{thm}

By Propositions \ref{Prop2.2} and  \ref{wc}, in case $f$ is weakly convex and $g$ is quasi-convex in variables $(x,y)$, we have $\partial^cv_\gamma(\bar x,\bar y;u,v)=\partial v_\gamma(\bar x,\bar y;u,v)$ for sufficiently small $\gamma$.
Hence by  Proposition \ref{wcsbdiff} and Theorem \ref{opt0}, we obtain the following optimality conditions for this case immediately. If 
in addition around an optimal solution $(\bar x,\bar y)$ of the bilevel program, the problems (VP) and (VP)$_\gamma$ are equivalent, then 
the following theorem also gives a necessary optimality condition for the original bilevel problem.

\begin{thm}\label{optwc}
Let $(\bar x,\bar y)$ be a local minimizer of {\rm (VP)$_\gamma$} for which $g(x,y)$ is quasi-convex in variables $(x,y)$ on $\mathbb R^n\times\mathbb R^m$, and $f$ is $\rho_f$-weakly convex on $\operatorname{gph}{\cal F}$ with $\rho_f>0$. 
  Suppose that the $\mathcal F(x)\neq\emptyset$ for all $x$ in a neighborhood of $\bar x$ and  MSCQ for system $g( x,y)\leq0$ holds at $(\bar x, \bar y)$. Then
$$
	 C(\bar x,\bar y)=\left \{(u,v) ~\Bigg| ~\begin{array}{l}
		 \nabla f(\bar x,\bar y)(u,v)=  \max_{\lambda\in \Lambda(\bar x,\bar y)} \nabla_x {\cal L}(\bar x,\bar y;\lambda)u,
	\quad 	\nabla F(\bar x,\bar y)(u,v)\leq0,\\
		\nabla g_i(\bar x,\bar y)(u,v)\leq 0,\ \forall i\in I_{g}(\bar x,\bar y), {\nabla G_i}(\bar x,\bar y)(u,v)\leq 0,\ \forall i\in I_{G}(\bar x,\bar y)\end{array} \right \}.
 $$
If {\rm (VP)$_\gamma$} for some $\gamma\in(0,1/\rho_f)$ is partially calm at $(\bar x,\bar y)$ in a direction $(u,v)\in C(\bar x, \bar y)$. 
Then there exists a vector $(\alpha,\lambda_g,\lambda_G)\in \mathbb{R}^{1+p+q}$ and $\bar \lambda\in\Lambda(\bar x,\bar y;u,v)$ satisfying 
		\begin{equation}\label{wcKKT}
		\begin{aligned}
			&0=\nabla  F(\bar x,\bar y)+\nabla g(\bar x,\bar y)^T (\lambda_g-\alpha \bar \lambda) +\nabla G(\bar x,\bar y)^T\lambda_G, \\
			&\alpha\geq0, 0\leq\lambda_g\perp g(\bar x,\bar y), \lambda_g\perp\nabla g(\bar x,\bar y)(u,v),
			0\leq\lambda_G\perp G(\bar x,\bar y), \lambda_G\perp\nabla G(\bar x,\bar y)(u,v).
		\end{aligned}
	\end{equation}
\end{thm}

Similarly, by Proposition \ref{ddv} under the RCR-regularity, one can obtain the following result. Note that RCR-regularity implies MSCQ  for the system $g( x,y)\leq0$.
\begin{thm}\label{optR}
	Let $(\bar x,\bar y)$ be a local minimizer of {\rm (VP)$_\gamma$}. 
	Suppose that the feasible map $\mathcal F$ satisfies RCR-regularity and RS at $(\bar x,\bar y)$ in a direction $u\in \mathbb{R}^n$. Then
	 $$ C(\bar x,\bar y)=\left \{(u,v) ~\Bigg| ~\begin{array}{l}
		 \nabla f(\bar x,\bar y)(u,v)=  \max_{\lambda\in \Lambda(\bar x,\bar y)} \nabla_x {\cal L}(\bar x,\bar y;\lambda)u,
	\quad 	\nabla F(\bar x,\bar y)(u,v)\leq0,\\
		\nabla g_i(\bar x,\bar y)(u,v)\leq 0,\ \forall i\in I_{g}(\bar x,\bar y), {\nabla G_i}(\bar x,\bar y)(u,v)\leq 0,\ \forall i\in I_{G}(\bar x,\bar y)\end{array} \right \}.
 $$
  Furthermore, suppose that {\rm (VP)$_\gamma$} for some $\gamma>0$ is partially calm at $(\bar x,\bar y)$ in direction $(u,v)\in C(\bar x,\bar y)$. Then there exists a vector $(\alpha,\lambda_g,\lambda_G)\in \mathbb{R}^{1+p+q}$ and $\bar \lambda\in\Lambda(\bar x,\bar y;u,v)$ satisfying (\ref{wcKKT}).
%
\end{thm}
{
\beginproof
The proof follows from Theorem \ref{opt0}, Proposition \ref{ddv} and the implication $$(u,v)\in C(\bar x,\bar y)\implies v\in\mathcal C(\bar  x,  \bar y;u),$$
where 
\begin{align*}
		& \mathcal C(\bar x, \bar y;u)= \left \{d ~\Bigg| \begin{array}{l}
			\max_{\lambda\in\Lambda(\bar x,\bar y,\bar w)}\nabla_x{\cal L}(\bar x,\bar y;\lambda)u	 = \nabla f(\bar x,\bar w)(u,d) \\ \nabla g_i(\bar x,\bar y)(u,d)\leq 0, \ \ \forall i\in I_g(\bar x,\bar y) \end{array} 
		\right \}.
	\end{align*} 
\endproof
}
By Corollary \ref{cor4.1} and Theorem \ref{opt0}, we obtain the following optimality conditions for this case immediately.  Note that MFCQ for the system $g(\bar x,y)\leq0$ at $\bar y$ implies MSCQ  for the system $g( x,y)\leq0$ at $(\bar x,\bar y)$.
\begin{thm}\label{optm}
Let $(\bar x,\bar y)$ be a local minimizer of {\rm (VP)$_\gamma$}.
	Suppose that either MFCQ for the system $g(\bar x, y)$ holds at $y=\bar y$ or condition (a) in Corollary \ref{cor4.1} holds. 
 Furthermore, suppose that {\rm (VP)$_\gamma$} for some $\gamma>0$ is partially calm at $(\bar x,\bar y)$ in direction $(u,v)\in C(\bar x,\bar y)$. Then there exists a vector $(\alpha,\lambda_g,\lambda_G)\in \mathbb{R}^{1+p+q}$ { together with $\bar \lambda\in\Lambda(\bar x,\bar y)$ satisfying (\ref{wcKKT}).}
 	\end{thm}
\beginproof
By Corollary \ref{cor4.1}, 
(\ref{upperestimate}) holds.
Meanwhile, for any direction $(u,d)$, $W(\bar x,\bar y,\bar y;u,d)\subseteq W(\bar x,\bar y,\bar y;0,0)$, which is convex. This together with (\ref{upperestimate}) and Definition \ref{Clarke} imply that $\partial^cv_\gamma(\bar x,\bar y;0,0)\subseteq W(\bar x,\bar y,\bar w;0,0)$. Then the proof is complete noticing that $\Lambda(\bar x,\bar y;0,0)=\Lambda(\bar x,\bar y)$.
\endproof

In particular if we take $(u,v)=(0,0)$ in (\ref{wcKKT}), then by Theorems \ref{optwc}, \ref{optR} and \ref{optm} we have the following corollary.
\begin{cor}\label{cor5.1}
Let $(\bar x,\bar y)$ be a local minimizer of {\rm (VP)$_\gamma$}.
	Suppose one of the following conditions holds:
	\begin{itemize}
	\item[(i)] $g(x,y)$ is quasi-convex, $f$ is $\rho_f$-weakly convex with $\rho_f>0$ on $\operatorname{gph}{\cal F}$. ${\cal F}(x) \not =\emptyset$ for all $x$ in a neighborhood of $\bar x$. MSCQ holds for the system $g(x,y)\leq 0 $ at $(\bar x,\bar y)$.
	\item[(ii)]  the feasible map $\mathcal F$ satisfies RCR-regularity and RS at $(\bar x,\bar y)$.
	\item[(iii)] MFCQ holds for the system $g(\bar x,y)\leq 0 $ at $\bar y$. 
	\end{itemize} Suppose that {\rm (VP)$_\gamma$}  for some $\gamma>0$ {($\gamma\in(0,1/\rho_f)$ when condition (i) holds)} is partially calm at $(\bar x,\bar y)$ in certain direction $(u,v)\in  C(\bar x,\bar y)$. Then there exists a vector $(\alpha,\lambda_g,\lambda_G)\in \mathbb{R}_+\times \mathbb{R}^{p+q}$ and $\bar \lambda\in\Lambda(\bar x,\bar y)$  satisfying 
	\begin{eqnarray}\label{sKKT}
		\begin{aligned}
			&0=\nabla  F(\bar x,\bar y)+\nabla g(\bar x,\bar y)^T (\lambda_g-\alpha \bar \lambda) +\nabla G(\bar x,\bar y)^T\lambda_G,\\
			&0\leq\lambda_g\perp g(\bar x,\bar y),\quad  0\leq\lambda_G\perp G(\bar x,\bar y).\
		\end{aligned}
	\end{eqnarray}
\end{cor}
In Theorem \ref{optm} and Corollary \ref{cor5.1}(iii), we need  a direction $(u,v)\in C(\bar x,\bar y)$. Such a direction can be obtained for example if the multiplier set $\Lambda(\bar x,\bar y)=\{\bar \lambda \} $ is a singleton since in this case the critical cone  is equal to
$$
	 C(\bar x,\bar y)=\left \{(u,v) ~\Bigg| ~\begin{array}{l}
		 \nabla f(\bar x,\bar y)(u,v)=   \nabla_x {\cal L}(\bar x,\bar y;\bar \lambda)u,
	\quad 	\nabla F(\bar x,\bar y)(u,v)\leq0,\\
		\nabla g_i(\bar x,\bar y)(u,v)\leq 0,\ \forall i\in I_{g}(\bar x,\bar y), {\nabla G_i}(\bar x,\bar y)(u,v)\leq 0,\ \forall i\in I_{G}(\bar x,\bar y)\end{array} \right \}.
 $$
\begin{remark} 	 The {stationarity} condition    (\ref{sKKT})  is in fact the {stationarity} condition (\ref{eqn2}).  Hence Corollary \ref{cor5.1} has provided sufficient conditions under which  the 
{stationarity} system (\ref{eqn2}) is a necessary optimality condition. 
{In fact the {stationarity} system (\ref{eqn2})  was proved  to be a necessary optimality condition for the class of jointly convex lower-level program under the nondifferentiable Abadie CQ in \cite[Theorem 4.2]{Ye04}, the extended Abadie CQ in \cite[Theorem 4.1 and Corollary 4.1]{Ye06} and  the  calmness of the partially perturbed feasible region of  $({VP})$
 (see \cite[Lemma 5.1(c)]{LM2023b}). }
 The directional version  (\ref{wcKKT}) is sharper than the nondirectional one  (\ref{sKKT}), and the directional S-{stationarity} conditions for the MPCC approach, see e.g., \cite[Definition 4.4]{Liang}. Note that unlike the KKT approach where $\bar \lambda$ in Definition \ref{SMPCC} is arbitrary, our optimality conditions come with a  lower-level multiplier $\bar \lambda$ which is selected and not given a priori.  
\end{remark}


 It is known that the directional quasi-normality and the relaxed constant positively linear dependence (RCPLD) condition (see \cite{XY}) can hold at feasible points of (VP).  
 Since RCPLD or the quasi-normality implies MSCQ, according to the discussion in Section \ref{Sec3}, the  MSCQ of the constraint system for problem {\rm (VP)$_\gamma$} in a direction  implies the  calmness of the problem in the same direction. Consequently, the directional partial calmness assumed in Theorems \ref{optwc}, \ref{optR} and \ref{optm} can be replaced by the directional quasi-normality condition, or the RCPLD condition. 
 
In the rest of this section, we will propose some verifiable sufficient conditions for the directional quasi-normality.   For the sake of simplicity, we only give the definition of directional quasi-normality condition for (VP)$_\gamma$. Recall that when  the lower-level problem is  generalized convex, 
(VP) is always equivalent to ${\rm (VP)}_\gamma$. However even without the generalized convexity assumption,  since $v(x)\leq v_\gamma(x,y)$ is always true, an optimal solution  $(\bar x,\bar y)$ of problem (VP) 
	must be  a feasible solution of problem ${\rm (VP)}_\gamma$. Hence we propose the following definition for any feasible solution  of {\rm (VP)$_\gamma$}. 
\begin{defn}\label{Defn5.2}
Let $(\bar x,\bar y)$ be a feasible solution  of {\rm (VP)$_\gamma$}. Assume that $v_\gamma(x,y)$ is directionally Lipschitz continuous at $(\bar x,\bar y)$ in direction $(u,v)$.
	We say that the directional quasi-normality for problem {\rm (VP)$_\gamma$} holds at $(\bar x,\bar y)$ in direction $(u,v)$ 
	 if there exists no nonzero vector $(\alpha,\nu_g,\nu_G)\in\mathbb R^{1+p+q}_+$ satisfying that \begin{equation}\label{KKT1}
	\begin{aligned}
				&0\in \alpha \nabla f(\bar x,\bar y)-\alpha \partial^c v_\gamma(\bar x,\bar y;u,v)+\nabla g(\bar x,\bar y)^T \nu_g
			+\nabla G(\bar x,\bar y)^T \nu_G,\\
			&\nu_g\perp g(\bar x,\bar y),\ \nu_g\perp\nabla g(\bar x,\bar y)(u,v),\  \nu_G\perp G(\bar x,\bar y),\ \nu_G \perp \nabla G(\bar x,\bar y)(u,v), \\
			 &\alpha \perp \xi \mbox{ for some } \xi \in \nabla f(\bar x,\bar y)(u,v)-D v_\gamma(\bar x,\bar y)(u,v) 
	\end{aligned}
	\end{equation}
	and there exist sequences $t_k\downarrow0,\ (u^k,v^k)\rightarrow(u,v)$ such that 
	\begin{eqnarray}\label{sequencial3}
		\begin{aligned}
		\alpha(f(\bar x+t_ku^k,\bar y+t_kv^k)-v_\gamma(\bar x+t_ku^k,\bar y+t_kv^k))>0,\ &\mbox{if}\ \alpha>0,\\
		g_i(\bar x+t_ku^k,\bar y+t_kv^k)>0,\ &\mbox{if}\ (\nu_g)_i> 0,  i\in I_g(\bar x,\bar y),  \\
		G_i(\bar x+t_ku^k,\bar y+t_kv^k)>0,\ &\mbox{if}\ (\nu_G)_i> 0, i\in I_G(\bar x,\bar y).
		\end{aligned}
	\end{eqnarray} 
\end{defn}
\begin{remark} 
Define $\phi(x,y):=(f(x,y)-v_\gamma(x,y),g(x,y),G(x,y))$ and $\lambda:=(\alpha,\nu_g,\nu_G)$. We can show that  the directional quasi-normality defined in Definition \ref{Defn5.2} is equivalent to the one defined in {Definition} \ref{qp} for the system $\phi(x,y)\leq 0$. By assumption, $\phi (x,y)$ is directionally Lipschitz continuous at $(\bar x,\bar y)$ in direction $(u,v)$.  
Then since $f(\bar x,\bar y)-v_\gamma(\bar x,\bar y)=0, g(\bar x,\bar y)\leq0, G(\bar x,\bar y)\leq0$, $(\ref{KKT1})$ means $0\leq \lambda\perp \phi(\bar x,\bar y)$ and $$\lambda\perp (\xi,\nabla g(\bar x,\bar y)(u,v), \nabla G(\bar x,\bar y)(u,v)),$$ 
where $(\xi,\nabla g(\bar x,\bar y)(u,v), \nabla G(\bar x,\bar y)(u,v))\in D\phi(\bar x,\bar y)(u,v). $ Since  by the calculus rule in \cite[Theorem 5.6]{Long}, we have
\begin{eqnarray*}
	\lefteqn{\partial (f-v_\gamma)(\bar x,\bar y; u, v)=\nabla f(\bar x,\bar y)+\partial(-v_\gamma)(\bar x,\bar y;u,v) }   \\
	&& \subseteq \nabla f(\bar x,\bar y)+\partial^c(-v_\gamma)(\bar x,\bar y;u,v)  = \nabla f(\bar x,\bar y)-\partial^c v_\gamma(\bar x,\bar y;u,v).
\end{eqnarray*} 
we conclude that (\ref{KKT1})-(\ref{sequencial3}) is stronger than the directional quasi-normality defined in Definition \ref{qp}. \end{remark}


Applying  the formulas of the directional derivatives  and the upper estimates for the directional subdifferentials for the function $v_\gamma(x,y)$ obtained in Section 4, we  derive  several sufficient conditions for the directional quasi-normality  below.
\begin{thm}\label{opt1}	Let $(\bar x,\bar y)$ be a feasible solution of {\rm (VP)$_\gamma$} and $u\in \mathbb{R}^n$.

	\begin{itemize}
		\item[(i)]
		Suppose either MFCQ holds for the system $g(\bar x,y)\leq0$ at $\bar y$, or the following conditions hold:
			\begin{itemize}
				\item The system $g(x,y)\leq0$ satisfies MSCQ at  $(\bar x,\bar y)$ and $\mathcal F$ satisfies RS at $(\bar x,\bar y)$ in direction $u$. $\mathbb L(\bar x,\bar y;\pm u)$ is nonempty and FOSCMS holds at $\bar y$ in direction $\nabla g(\bar x,\bar y)(\pm u,v)$ for any $v\in\mathbb L(\bar x,\bar y;\pm u)$. 
	
			\end{itemize}	Suppose there exists no nonzero vector $(\alpha,\nu_g,\nu_G)\in\mathbb R^{1+p+q}_+$ satisfying conditions (\ref{KKT1})-(\ref{sequencial3}) with $\partial^c v_\gamma(\bar x,\bar y;u,v)$ replaced by 
\[
\operatorname{co}\left( \bigcup_{\lambda\in \Lambda(\bar x,\bar y;u,d) \atop d\in\widetilde {\mathcal C}(\bar x,\bar y; u)\cup(\widetilde {\mathcal C}(\bar x, \bar y;0)\cap\mathbb S)} 
\left\{(\nabla_x f(\bar x,\bar y)+\nabla_x g(\bar x,\bar y)^T\lambda)\right\} \times \{0\}\right),
\]
where
\begin{align*}
& \widetilde {\mathcal C}(\bar x,\bar y; u)\\
:=&\left \{d ~\Bigg| \begin{array}{l}
\min_{\lambda\in \Lambda(\bar x,\bar y)} \nabla_x 
{\cal L}(\bar x,\bar y;\lambda)u  \leq \nabla f(\bar x,\bar y)(u,v)\leq  \max_{\lambda\in \Lambda(\bar x,\bar y)} \nabla_x {\cal L}(\bar x,\bar y;\lambda)u\\\nabla g_i(\bar x,\bar y)(u,d)\leq 0, \ \ \forall i\in I_g(\bar x,\bar y) \end{array} 
\right \}.
\end{align*}
Moreover if $S_\gamma(x,y)$ is inner calm* at $(\bar x,\bar y)$ in direction $(u,v)$ then $\partial^c v_\gamma(\bar x,\bar y;u,v)$ is replaced by \[
\operatorname{co}\left(\bigcup_{\lambda\in \Lambda(\bar x,\bar y;u,d)\atop d\in\widetilde {\mathcal C}(\bar x,\bar y; u)} \left\{(\nabla_x f(\bar x,\bar y)+\nabla_x g(\bar x,\bar y)^T\lambda)\right\} \times \{0\}\right).
\] Then the directional quasi-normality for problem {\rm (VP)$_\gamma$} holds at $(\bar x,\bar y)$ in direction $(u,v)$.
	\item[(ii)] Suppose $f$ is $\rho_f$-weakly convex on $\operatorname{gph}{\cal F}$ with $\rho_f>0$. Suppose that  $\mathcal F(x)\neq\emptyset$ for all $x$ in a neighborhood of $\bar x$ and MSCQ for the system $g( x,y)\leq 0$ holds at $(\bar x, \bar y)$.  Suppose there exists no nonzero vector $(\alpha,\nu_g,\nu_G)\in\mathbb R^{1+p+q}_+$ satisfying conditions (\ref{KKT1})-(\ref{sequencial3}) with $
	\partial^c v_\gamma(\bar x,\bar y;u,v)$ replaced by $$\bigcup_{\lambda_g \in\Lambda(\bar x,\bar y;u,v)} (\nabla_xf(\bar x,\bar y)+\nabla_x g(\bar x,\bar y)^T\lambda_g,0).$$ 	Then the directional quasi-normality for problem {\rm (VP)$_\gamma$} holds at $(\bar x,\bar y)$ in direction $(u,v)$.	
	\item[(iii)] Suppose that $\bar y\in S(\bar x)$ and
	${\mathcal F}$ is RCR-regular and satisfies RS at $(\bar x,\bar y)$ in direction $u$.
	 Suppose there exists no nonzero vector $(\alpha,\nu_g,\nu_G)\in\mathbb R^{1+p+q}_+$ satisfying conditions (\ref{KKT1})-(\ref{sequencial3}) with $\partial^cv_\gamma(\bar x,\bar y;u,v)$ replaced by  $$\bigcup_{\lambda_g \in\Lambda(\bar x,\bar y;u,d)} (\nabla_xf(\bar x,\bar y)+\nabla_x g(\bar x,\bar y)^T\lambda_g,0),$$
	 where  $d\in \mathcal C(\bar x,\bar y; u)$.
   Then the directional quasi-normality for problem {\rm (VP)$_\gamma$} holds at $(\bar x,\bar y)$ in direction $(u,v)$.
	\end{itemize}
\end{thm}        
\beginproof
By Proposition \ref{dini}, $\widetilde C(\bar x,\bar y;u)$ contains $C(\bar x,\bar y;u)$ in case (i). The rest of the proof follows from Definition \ref{Defn5.2} and upper estimates for $\partial^cv_\gamma(\bar x,\bar y;u,v)$ in  Propositions  \ref{wcsbdiff}, \ref{est1} and \ref{ddv}.
\endproof

In the following example, we verify the directional quasi-normality using Theorem \ref{opt1} (ii). {Note that although the lower-level generalized convexity fails, the global minimizer $(\bar x,\bar y)$ of (VP) remains a local minimizer of (VP)$_\gamma$, hence optimality conditions of (VP)$_\gamma$ still characterize $(\bar x,\bar y)$. Besides, we can verify the directional quasi-normality for problem {\rm (VP)$_\gamma$}  at $(\bar x,\bar y)$.}

\begin{example} Consider the following bilevel program
	\begin{eqnarray*}
		\quad & \displaystyle \min_{x,y}  & F(x,y):=(x-y)^2\\
		& {\rm s.t.}&  y \in S(x):=\arg\min_y \left\{-\left(x-y\right)^2\left| 
		\begin{array}{ll}
			y-x-1\leq0,\\
			x-y-1\leq0
		\end{array}\right.\right\}.
	\end{eqnarray*}
	Denote by $f(x,y):=-(x-y)^2, g_1(x,y):=y-x-1, g_2(x,y):=x-y-1.$	It is easy to verify 
	\begin{equation}\label{esolution}
		S(x)=\{x+1, x-1\} \mbox{ and }
		v(x)=-1.
	\end{equation} 
	Point $(\bar x,\bar y)=(0,-1)$ is a global optimal solution of the bilevel program.

	
	\begin{equation*}
		\nabla F(\bar x,\bar y)=\left [\begin{matrix}
			2\\
			-2
		\end{matrix}\right ], \quad \nabla f(\bar x,\bar y)=\left [\begin{matrix}
			-2\\
			2
		\end{matrix}\right ] \quad \nabla g_1(\bar x,\bar y)=\left [\begin{matrix}
			-1\\
			1
		\end{matrix}\right ],  \quad \nabla g_2(\bar x,\bar y)=\left [\begin{matrix}
			1\\
			-1
		\end{matrix}\right ].
	\end{equation*}
	Note that $f$ is weakly convex with $\rho_f=2$, $g_i(i=1,2)$ is convex, and $\mathcal F(x)\neq\emptyset$ near $\bar x$.	Since $g(x,y)$ is linear, MSCQ for the system $g(x,y)\leq0$ holds at $(\bar x,\bar y)$.
	Since $g_1(\bar x,\bar y)=-2<0, g_2(\bar x,\bar y)=0$, the set of multipliers and the set of directional multipliers in direction $(u,v)$ are  
	\begin{eqnarray}
		\Lambda(\bar x,\bar y;u,v)=\Lambda(\bar x,\bar y)=\{(0,2)\}.
	\end{eqnarray}
	Then taking $\gamma\in(0,1/2)$, the feasible region of {\rm (VP)}$_\gamma$ is set $$\{(x,y)|y=x+1\mbox{ or } y=x\mbox{ or } y=x-1\}.$$
{One can easily verify, in a small neighborhood of $(\bar x,\bar y)$, say $\mathbb B_{r}(\bar x,\bar y)$ with $0<r<0.1$, the set of feasible points of {\rm (VP)} and {\rm (VP)}$_\gamma$ are both the set $\{(x,y)\in\mathbb B_{r}(\bar x,\bar y)|y=x-1\}$. Then $(\bar x,\bar y)$ is a local minimizer of {\rm (VP)}$_\gamma$. We know the global solution of {\rm (VP)}$_\gamma$ is the set $\{(x,y)|y=x\}$, hence $(\bar x,\bar y)$ is only a local minimizer.}
	
By {Propositions \ref{wc} and \ref{wcsbdiff}}, $v_\gamma(x,y)$ is directionally Lipschitz continuous and directionally differentiable in direction $(u,0)$ with $u>0$ and 
	\begin{align*}
		v'_\gamma(\bar x,\bar y;u,0)
		=&\max_{\lambda\in{\Lambda(\bar x,\bar y;u,0)}}\nabla_x\mathcal L(\bar x,\bar y,\bar y;\lambda)u=0.
	\end{align*}
	Now we prove that the directional quasi-normality for problem {\rm (VP)$_\gamma$} holds at $(\bar x,\bar y)$.
	The critical cone can be calculated as
	\begin{align*}
		C(\bar x,\bar y)&:=\{(u,v)|\nabla F(\bar x,\bar y)(u,v)\leq0, \nabla f(\bar x,\bar y)(u,v)-v'_\gamma(\bar x,\bar y;u,v)=0,\nabla g_2(\bar x,\bar y)(u,v)\leq 0\} \\
		&= \{(u,v)| u-v\leq0, \nabla f(\bar x,\bar y)(u,v)-v'_\gamma(\bar x,\bar y;u,v)=0\}\\
		&= \{(u,v)| u-v\leq0, \nabla f(\bar x,\bar y)(u,v)-v'_\gamma(\bar x,\bar y;u,0)=0\}\\
		&= \{(u,v)| u-v\leq0, -2u+2v=0\}\\
		&= \{(u,v)| u-v= 0\},
	\end{align*} {where the third equality follows from the fact $ S_\gamma(\bar x,\bar y)=\{\bar y\}$ and (\ref{ddformula}).}
	Let $(\bar u,\bar v)=(1,1)$, we have  $(\bar u,\bar v)\in C(\bar x,\bar y)$. 
	Since $g_1(\bar x,\bar y)<0=g_2(\bar x,\bar y), \nabla g_2(\bar x,\bar y)(\bar u,\bar v)=0$, by Proposition \ref{wcsbdiff} we have
	\begin{eqnarray*}
		\partial_x^cv_\gamma(\bar x,\bar y;\bar u,\bar v)=  \left \{\nabla_x f(\bar x,\bar y)+\nabla_x g(\bar x,\bar y)^T\lambda_g| \lambda_g \in \Lambda(\bar x,\bar y;\bar u,\bar v) \right\}=0
	\end{eqnarray*} and $\partial_y^cv_\gamma(\bar x,\bar y;\bar u,\bar v)=\{0\}$.
	Let 
	$\alpha, \nu
	$ be such that
	\begin{eqnarray}\label{e3x}
		&& 0\in\alpha(\nabla_x f(\bar x,\bar y)-\partial_x^cv_\gamma(\bar x,\bar y;\bar u,\bar v))+\nu\nabla_xg_2(\bar x,\bar y),\label{e3z} \\
		&&	0=\alpha\nabla_y f(\bar x,\bar y)+\nu\nabla_yg_2(\bar x,\bar y),\label{e3y}\\
		&& {\nu\nabla g_2(\bar x,\bar y)(\bar u,\bar v)=0}, \alpha \geq0, \nu\geq0 
	\end{eqnarray}
	and there exist sequences $t_k\downarrow0,\ (u^k,v^k)\rightarrow (\bar u,\bar v)$, such that
	\begin{eqnarray}
		&&f(\bar x+t_ku^k,\bar y+t_kv^k)-v_\gamma(\bar x+t_ku^k,\bar y+t_kv^k)>0\ \mbox{if}\ \alpha>0,\label{eseq1}\\
		&& g_2(\bar x+t_ku^k,\bar y+t_kv^k)>0\ \mbox{if}\ \nu>0.\label{eseq5}
	\end{eqnarray}
	(\ref{e3z}) and (\ref{e3y}) imply that $-2\alpha+\nu=0$. Hence any $(\alpha, \nu)$ with $2\alpha=\nu$ satisfies (\ref{e3z})-(\ref{e3y}).  
	We now show the conditions (\ref{e3x})-(\ref{eseq5}) can only hold if $\alpha=\nu=0$.  To the contrary, assume $\alpha>0$. Since $\nu=2\alpha>0$, (\ref{eseq5}) is required to hold.
	Let  $t_k\downarrow0,\ (u^k,v^k)\rightarrow (\bar u,\bar v)$ be arbitrary and suppose that (\ref{eseq5}) holds. Then $g_2(x^k,y^k)>0$ for $(x^k,y^k):=(\bar x+t_ku^k,\bar y+t_kv^k)$. It follows that   $y^k<x^k-1$. Since $\nabla_y f(x^k,y^k)=2(x^k-y^k)>2$ and $y^k<x^k-1$, we have  {$f(x^k,y^k)<f(x^k,x^k-1)=v(x^k)\leq v_\gamma(x^k,y^k)$,} where the last equality follows from  $(\ref{esolution})$. Hence (\ref{eseq1}) does not hold.
	The contradiction shows that  $(\alpha,\nu)=(0,0)$ and following Theorem \ref{opt1} (ii), the directional quasi-normality  for problem {\rm (VP)$_\gamma$} holds at $(\bar x,\bar y)$ in direction $(\bar u,\bar v)$.
\end{example}

\section{Conclusion}

The classical approach of solving a bilevel program by its MPCC reformulation is known to be problematic when the lower-level program has non-unique multipliers.  In particular a local optimal solution $(\bar x,\bar y, \bar \lambda)$ for MPCC may not imply that $(\bar x,\bar y)$ is a local optimal solution of the bilevel program. Recently using the value function reformulation, some optimality conditions such as system (\ref{eqn2}) which does not depend on an arbitrarily chosen multiplier of the lower-level program, have been derived and used as a basis for numerical computation. However the validity of system (\ref{eqn2}) as an optimality condition for bilevel program requires very strong assumptions such as the joint convexity of the function $f$ and $g$, as well as others. In this paper by using Moreau envelope reformulation, we are able to obtain sufficient conditions under which system (\ref{eqn2}) is an optimality condition.

\section{Acknowledgments}

The authors are grateful to the associate editor and two referees for their helpful comments
and constructive suggestions.

\end{document}